\documentclass{amsart}

\usepackage{amsfonts,amssymb,amsmath,amsthm}
\usepackage{cite}
\usepackage{dptsymb}
\usepackage[dvips, final]{graphics}
\usepackage{epsfig}
\usepackage{subfigure}
\usepackage[latin1]{inputenc}
\usepackage[T1]{fontenc}
\usepackage{url}

\def\printname#1{
        \if\draft y
                \smash{\makebox[0pt]{\hspace{-0.5in}
                        \raisebox{8pt}{\tt\tiny #1}}}
        \fi
}

\newcommand{\eepic}[2]{}
\newcommand{\silenteepic}[2]{}

\newread\testin
\def\maybeinput#1{
\openin\testin=#1
\ifeof\testin\typeout{Warning: input #1 not found}\else\input#1\fi
\closein\testin
}

\newwrite\transout
\openout\transout=figs.list
\newcommand{\figscaled}[2]{%
\write\transout{-m #2 draws/#1.fig}%
\maybeinput{draws/#1.tex}%
}

\newcommand{\figcent}[2]{\mathcenter{\figscaled{#1}{#2}}}

\usepackage{calc}
\def\mathcenter#1{
        \raisebox{.5ex+\depth-.5\totalheight}{\hbox{#1}}
}

\newtheorem{Definition}{Definition}
\newtheorem{Lemma}{Lemma}
\newtheorem{Remark}{Remark}
\newtheorem{Theorem}{Theorem}
\newtheorem{Proposition}{Proposition}

\theoremstyle{definition}
\newtheorem*{Ack}{Acknowledgment}
\newcommand{\R}{\mathbb R}

\newcommand{\N}{\mathbb N}
\newcommand{\C}{\mathbb C}

\newcommand{\op}{{\mathcal O}}

\newcommand{\oploc}[1]{\op_{\mathrm{loc}}^{#1}}
\newcommand{\opprod}[1]{\op_{\Delta}^{#1}}
\newcommand{\opform}[1]{\op_{\mathrm{form}}^{#1}}
\newcommand{\opdef}[1]{\op_{\mathrm{def}}^{#1}}

\newcommand{\vect}{\operatorname{Vect}}

\newcommand{\cotR}[1]{\mathrm{T}^*\R^{#1}}
\newcommand{\locfuncRd}[1]{\mathcal O^{#1}_{\mathrm{loc}}(\cotR d)}

\newcommand{\id}{\operatorname{id}}
\newcommand{\morph}{\operatorname{Hom}}
\newcommand{\obj}{\operatorname{Obj}}
\newcommand{\rtree}[1]{\figcent{#1}{0.4}}

\newcommand{\graph}{\operatorname{graph}}

\newcommand{\neutral}{\textbf{e}}

\newcommand{\ii}{{\mathrm{i}}}
\newcommand{\dd}{{\mathrm{d}}} 
\newcommand{\ee}{{\mathrm{e}}}

\usepackage{verbatim}

\title{\bf Formal Lagrangian Operad}

\author[A.~S.~Cattaneo]{Alberto~S.~Cattaneo}
\address{Institut f\"ur Mathematik, Universit\"at Z\"urich--Irchel, Winterthurerstrasse 190, CH-8057 Z\"urich, Switzerland}
\email{alberto.cattaneo@math.unizh.ch}

\author[B.~Dherin]{Benoit Dherin}
\address{D-MATH, ETH-Zentrum, CH-8092 Z\"urich, Switzerland}
\email{dherin@math.ethz.ch}

\author[G.~Felder]{Giovanni Felder}
\address{D-MATH, ETH-Zentrum, CH-8092 Z\"urich, Switzerland}
\email{felder@math.ethz.ch}

\thanks{A.~S.~C. acknowledges partial support of SNF Grant No.~200020-107444/1 and the IHES for hospitality}
\thanks{B. D. and G. F. acknowledge partial support of SNF Grant No.~21-65213.01}

\begin{document}

\maketitle

\begin{abstract}

Given a symplectic manifold $M$, we may define an operad structure on the
the spaces $\op^k$ of the Lagrangian submanifolds of $(\overline{M})^k\times M$ via
symplectic reduction. If $M$ is also a symplectic groupoid, then its multiplication
space is an associative product in this operad. Following this idea,
we provide a deformation theory for symplectic groupoids analog to the deformation theory of algebras.
It turns out that the semi-classical
part of Kontsevich's deformation of $C^\infty(\R^d)$ is a deformation of the
trivial symplectic groupoid structure of $T^*\R^d$.

\end{abstract}

\tableofcontents

\section{Introduction}

Symplectic groupoids, in the extended symplectic category, may be thought as the analog
of associative algebras in the category of vector spaces. For the latter, a deformation theory
exists and is well known. In this article, we will present a conceptual framework as well
as an explicit deformation of the trivial symplectic groupoid over $\R^d$. In fact, rephrased
appropriately, most constructions of the deformation theory of algebras can be extended
to symplectic groupoids, at least for the trivial one over $\R^d$. Our guide line will be the Kontsevich deformation
of the usual algebra of functions over $\R^d$, $\big(C^\infty(\R^d),\cdot\big)$. Namely, the usual point-wise
product of functions $S_0^2(f,g) = fg$ generates a suboperad, the product suboperad,
$\op_S^n=\big\{S_0^n \big\},$ of the endomorphism operad $\op$ of 
$C^\infty(\R^d)$, where $S_0^n$ is the $n$-multilinear map defined by 
$S_0^n(f_1,\dots,f_n) = f_1f_2\dots f_n.$ For each $n$ one may choose the vector subspace
$\opdef n\subset \op^n$
of $n$-multidifferential operators. The operad structure of $\op$ induces 
an operad structure on $\op_S+\opdef{}$, which in turns generates an operad structure on $\opdef{}$
which is, however, non-linear. Then, $\gamma$ is a deformation of the usual product $S_0^2$, i.e.,  an element
$\gamma\in\opdef2$ such that $S_0^2+\gamma$ is still an associative product, iff $\gamma$
is a product in the induced deformation operad $\opdef{}$. We may also consider the formal
version by replacing $\opdef{}$ by the formal power series in $\epsilon$, $\epsilon \opdef{}[[\epsilon]]$.
M. Kontsevich in \cite{kontsevich1997} gives an explicit formal deformation of the product of
functions over $\R^d$,
$$S_\epsilon = S_0^2+\sum_{n=1}^\infty\epsilon^n\sum_{\Gamma\in G_{n,2}}W_\Gamma B_\Gamma,$$
where the $W_\Gamma$'s are the Kontsevich weights and the $B_\Gamma$'s the Kontsevich
bidifferential operators associated to the Kontsevich graphs of type $(n,2)$ (see 
\cite{CDF2005}
for a brief introduction
).

If we consider the trivial symplectic groupoid $\cotR d$ over $\R^d$, we see that
the multiplication space $$\Delta_2^n:=\Big\{(p_1,x),(p_2,x),(p_1+p_2,x):p_1,p_2\in \R^{d*},
x\in \R^d\Big\}$$ generates an operad $\opprod n = \big\{\Delta_n\big\},$
where $$\Delta_n := \Big\{(p_1,x),\dots,(p_n,x),(p_1+\dots+p_n,x):p_i\in\R^{d*},x\in\R^d\Big\}.$$
$\Delta_2$ is a product in this operad. The compositions are given by symplectic reduction as
the $\Delta_n$'s are Lagrangian submanifolds of $\overline{(\cotR d)^n}\times \cotR d$.
The main difference with the vector space case is that there is no ``true'' endomorphism
operad where $\opprod{}$ would naturally embed into. Thus, the question of finding
a deformation operad for $\opprod{}$ must be taken with more care. The first remark is that the
$\Delta_n$ may be expressed in terms of generating functions $$S_0^n(p_1,\dots,p_n,x) = (p_1+\dots+p_n)x.$$
Namely, $\Delta_n = \graph dS_0^n$. The idea is to look at the operad structure induced on the generating
functions by symplectic reduction. In fact it is possible to find a vector space of special functions 
$\opdef n$ for each $n$ such that $\opprod{}+\opdef{}$ remains an operad. The formal version of it
gives a surprising result. Namely, we may find an explicit deformation of the trivial generating function $S_0^2$
, it is given by the formula 
$$S_\epsilon = S_0^2 +\sum_{n=1}^\infty \epsilon^n\sum_{\Gamma\in T_{n,2}}W_\Gamma \hat B_\Gamma,$$
where the $W_\Gamma$ are the Kontsevich weights and the $\hat B_\Gamma$ are the symbols of the 
Kontsevich bidifferential operators and the sum is taken over all Kontsevich trees $T_{n,2}$. This formula
may be seen as the semi-classical part of Kontsevich deformation quantization formula.

As a last comment, note that Kontsevich derives its star product formula from a more general
result. In fact, he shows that $U = \sum_{n} \epsilon^n U_n$ where 
$$U_n(\xi_1,\dots,\xi_n) = \sum_{\Gamma\in G_n} W_\Gamma B_\Gamma(\xi_1,\dots,\xi_n)$$
for $\xi_i\in\Gamma(\wedge^{d^i}TM)$, $i=1,\dots,d$ 
is an $L_\infty$-morphism from the multivector fields to the multidifferential operators on $\R^d$.
In our perspective, we may still write 
$$\tilde U_n(\xi_1,\dots,\xi_n) = \sum_{\Gamma\in T_n} W_\Gamma  \hat B_\Gamma(\xi_1,\dots,\xi_n)$$
summing over Kontsevich trees instead of Kontsevich graphs and replacing multidifferential operators
by their symbols. Exactly, as in Kontsevich case, 
$$S_\epsilon = S_0^2 +\sum_{n\geq 1} \epsilon^n \tilde U_n(\alpha,\dots,\alpha)$$
is an associative deformation of the generating function of the trivial symplectic groupoid $T^*\R^d$.
However, it is still not completely clear how to define  ``semi-classical $L_\infty$-morphisms''.

\subsubsection*{Organization of the article}

In Section \ref{prodinext}, we describe the endomorphism operad $\op(M) = \morph(M^{\otimes^n},M)$
associated to any object $M$ in a monoidal category. We explain what is an associative product $S$ on $M$  in an monoidal category and 
we define the product suboperad $\op_S(M)$ of $\op(M)$.
If the category is further associative, we may choose a deformation operad for $S$, which is a choice, for
each $n\in \N$ of a vector subspace $\opdef n$ such that $\op_S+\opdef{}$ is still an operad. We describe the
deformations of $S$ in terms of products in $\opdef{}$. As an example of this construction, we expose
Kontsevich product deformation in this language. At last, we show that the extended symplectic category, although
not being a true category, exhibits monoidal properties allowing us to carry the precedent construction
up to a certain point. Then, we focus on the trivial symplectic groupoid over $\R^d$ case and
define the product operad associated to its multiplications space. We give a deformation operad
on a local form, the local deformation operad. In particular, we show that any local
deformation of the trivial product gives rise to a local symplectic groupoid over $\R^d$. We conclude
this Section by defining equivalence between deformations of the trivial generating function
and we show that two equivalent deformations induce the same local symplectic groupoid.

In Section \ref{combinatorics}, we describe the combinatorial tools needed to give a formal
version of the local Lagrangian operad. As the problem consists mainly in taking Taylor's series
of some implicit equations we need devices to keep track of all terms to all orders. The crucial
point is that these implicit equations, describing the composition in the local Lagrangian operad,
have a form extremely close to a special Runge-Kutta method: the partitioned implicit Euler method. 
We borrow then some techniques form numerical analysis of ODEs to make the expansion at all orders.

In the last Section, we describe the formal Lagrangian operad, which is the perturbative version of the
local one, in terms of composition of bipartite trees. We give in particular the product equation
in the formal deformation operad in terms of these trees. 
At last, we restate the main Theorem of \cite{CDF2005} in this language.  This tells us
that the semi-classical part of Kontsevich star product on $\R^d$ is a product in the formal deformation
operad of the cotangent Lagrangian operad in $d$ dimensions.

This article is inspired for a large part from unpublished notes \cite{Cattaneo2002} of one of the authors,
in which the notion of Lagrangian operad first appeared, and from the PhD thesis \cite{dherin2004} of an another author.
It is a natural development of results presented in \cite{CDF2005}.

\begin{Ack}We thank Domenico Fiorenza for useful comments and suggestions.\end{Ack}

\section{Product in  the extended symplectic category}\label{prodinext}

\subsection{Basic constructions and Kontsevich deformation}

In this Section, we describe, in any monoidal category, a natural generalization of an associative
algebra structure over a vector space. It is the notion of  product in the endomorphism operad $\op(M)$
of an object $M$ in the category. If the category is further additive, we explain what is a deformation of a
product $S\in\op^2(M)$ and construct a non-linear operad, the deformation operad $\op_{\mathrm{def}}(M,S)$ 
associated to $S$ in which any product is equivalent to a deformation of $S$. We present the well-known
Kontsevich deformation of the usual product of functions over $\R^d$ in this language. At last,
we see that most parts of this construction, can be applied to the extended symplectic category, leading
to the notion of Lagrangian operad.

\begin{Definition}
An operad $\op$ consists of 
\begin{enumerate}
\item a collection of sets $\op^n, n \geq 0$
\item composition laws
\begin{eqnarray*}
\op^n \times \op^{k_1} \times \dots \times \op^{k_n} 
& \longrightarrow &
 \op^{k_1 + \dots +k_n}\\
(F,G_1,\dots,G_n) & \mapsto & F(G_1,\dots,G_n)
\end{eqnarray*}
satisfying the following associativity relations,
\begin{gather*}
 F(G_1,\dots,G_n)(H_{11},\dots,H_{1k_1},\dots,H_{n1},\dots,H_{nk_n})
= \\
 F(G_1(H_{11},\dots,H_{1k_1}),\dots,G_n(H_{n1},\dots,H_{nk_n}))
\end{gather*}
\item a unit element $I \in \op^1$ such that 
$F(I,\dots,I) = F \textrm{ for all } F \in \op^n$
\end{enumerate}
It is usually also required some equivariant action of the symmetric group.
We do not require this here. 
\end{Definition}

The structure we have just defined should then be called more correctly ``non symmetric operad''. However, we will simply keep using
the term ``operad'' instead of ``non symmetric operad'' in the sequels.

\subsubsection*{Product in a monoidal category}

We consider here a monoidal category  $\mathcal{C}$. We denote by 
$\otimes : \mathcal{C} \times \mathcal{C} \longrightarrow \mathcal{C}$
the product bifunctor and  by $\neutral\in \mathcal{C}$ the neutral object.
Let us recall that we have the following canonical isomorphisms
$$(A \otimes B )\otimes C \simeq  A \otimes (B \otimes C)\quad\textrm{ and }\quad
\neutral \otimes A \simeq A \otimes \neutral \simeq A $$
for all $A,B,C\in\obj\mathcal C$.

Let $\mathcal{C}$ be a monoidal category and an object $M \in \obj\mathcal{C}$.
We define the endomorphism operad of $M$ in the following way:

\begin{enumerate}
\item $\op^n(M):=\morph(M^{\otimes n},M), \: \op^0(M) :=
\morph(\neutral ,M)$
\item $F(G_1,\dots,G_n):= F \circ (G_1 \otimes \dots \otimes G_n)$
\item the unit is given by $\id_M \in \op^1(M)$.
\end{enumerate}

The operad axioms follow directly from the bifunctoriality 
 of $\otimes$, i.e,
\begin{eqnarray*}
(f \otimes g) \circ (\psi \otimes \phi) & = & 
	(f \circ \psi)\otimes (g \circ \phi)\\
\id_M \otimes \dots \otimes \id_M & = & \id_{M \otimes \dots \otimes M}.
\end{eqnarray*}
If $M$ is an object of a monoidal category $\mathcal{C}$, we may define a 
product on $M$. 

\begin{Definition}
An associative product
\footnote{In \cite{GV1995}, Gerstenhaber and Voronov call it
a multiplication.}on an operad $\op$ is an element $S \in \op^2$ such that $S(I,S)=S(S,I)$.
An associative product on $M$ is an associative product in the endomorphism operad $\op(M)$.
In the sequels, we will constantly use the term product to mean in fact associative product.
\end{Definition}

Given a product $S\in\op^2$, the associativity of the operad implies that, for any 
$F \in \op^k$, $G \in \op^l$ and $H \in \op^m$ we have,
\begin{eqnarray*}
S(F,S(G,H)) & = & S (I,S)(F,G,H) \\
 & = & S(S,I)(F,G,H) \\
 & = & S(S(F,G),H).
\end{eqnarray*}

This notion is the natural generalization of an associative product on a vector
space. Namely, if $M$ is a vector space, $\op^2(M)$ is the set of bilinear map
 on $M$. As in this case $\op^0(M)=\morph(\C,M)=M$, we have that
$S : \op^0(M) \times  \op^0(M) \longrightarrow  \op^0(M)$ 
is an associative product on $M$. 

\subsubsection*{Product deformation in a monoidal additive category}

Suppose we have a product $S\in\op^2(M)$, where $M$ is an object of a monoidal 
category $\mathcal C$.  If the category $\mathcal{C}$ is further additive, 
we may try to deform $S$, i.e., to find an element 
$\gamma \in  \op^2(M)$ such that $S+\gamma$ is still a product.

At this point, we may follow the main trend and introduce the Hochschild complex 
of the linear operad $\op(M)$, define the bilinear Gerstenhaber bracket and
the Hochschild differential associated to the product $S$. A deformation of
$S$ would then be a solution of the Maurer-Cartan equation written in the Hochschild
Differential Graded Lie Algebra controlling the deformations of $S$. 

We will however rephrase slightly this deformation theory in a way allowing us
to include afterward the case when the morphisms of the category are still vector
spaces but the composition of morphism is no longer bilinear. 
That will be exactly the case we will have to deal with in the next Sections.

The first step is to notice that a product $S\in \op^2(M)$ generate a sub-operad
$\op_S(M)$, which we call a {\bf product operad}, in $\op(M)$ with only one point in each degrees:
\begin{gather*}
\op_S^0(M)  :=  \emptyset,\quad \op_S^1(M)  :=  \Big\{I\Big\},\quad \op_S^2(M)  :=   \Big\{S \Big\},\\
\op_S^3(M)  :=   \Big\{S(S,I) \Big\},\quad \op_S^4(M)  :=   \Big\{S(S(S,I),I) \Big\},\dots \textrm{etc}
\end{gather*}
To simplify the notation we will denote by $S_0^n$ the unique element in $\op_S^n(M)$.

\begin{Definition}
Let $M$ be an object of an additive  monoidal category $\mathcal C$ and
$S\in\op^2(M)$ a product. A \textbf{deformation operad}, $\op_{\mathrm{def}}(M,S)$ 
for $S$ is a choice for each $n\in\N$ of a vector space space 
$\opdef n(M,S)\subset\op^n(M),$
such that $\op_S+\op_{\mathrm{def}},$ is still an operad with respect to
the compositions of $\op(M)$, i.e., 
\begin{eqnarray}\label{operadeq}
  (S_0^n+\gamma)(S_0^{k_1}+\gamma^{1},\dots,S_0^{k_n}+\gamma^{n})  
  & = &   S_0^{k_1+\dots +k_n} +R(\gamma;\gamma^{1},\dots,\gamma^{n}),
\end{eqnarray}
with  $$R(\gamma;\gamma^{1},\dots,\gamma^{n})\in \opdef {k_1+\dots+k_n}(M,S),$$
for $\gamma\in\opdef n(M,S)$ and $\gamma^{i}\in\opdef {k_i}(M,S)$ for $i=1,\dots,n$.

We say that an element $\gamma\in\op_{\mathrm{def}}^2(M,S)$ is a {\bf deformation of the product} $S$ w.r.t.
the deformation operad $\op_{\mathrm{def}}$ if $S+\gamma$ is still a product in $\op_S+\op_{\mathrm{def}}$.
\end{Definition}

Remark that the $R$s are not multilinear.

\begin{Proposition}
Let $\op_{\mathrm{def}}(M,S)$ be a deformation operad for a product $S\in\op^2(M)$.
Then the compositions $$\gamma(\gamma^{1},\dots,\gamma^{n}):= R(\gamma;\gamma^{1},\dots,\gamma^{n}),$$
defined by equation (\ref{operadeq}) gives $\op_{\mathrm{def}}(M,S)$ together with the unit $0\in \op_{\mathrm{def}}^1(M,S)$ 
the structure of an operad.
\end{Proposition}

\begin{proof}
The proof is direct using only equation (\ref{operadeq}) and the operad structure of
the endomorphism operad $\op(M)$.
\end{proof}

\begin{Proposition}
Let $S\in \op^2(M)$ be a product. Take an element $\gamma\in \opdef 2(M,S)$. Then,
$\gamma$ is a deformation of the product $S$ iff $\gamma$ is a product in $\op_{\mathrm{def}}(M,S)$.
In particular, $0\in \opdef 2(M,S)$ is always a product in the deformation operad of $S$.
\end{Proposition}

\begin{proof}
$\gamma$ is a deformation of $S$ iff $$(S+\gamma)(S+\gamma,I) = (S+\gamma)(I,S+\gamma),$$
which is equivalent to $$S_0^3+R(\gamma;\gamma, 0) = S_0^3+R(\gamma;0,\gamma).$$
\end{proof}
From now on, we will write $0_1$ for the identity element of the deformation operad which is
the zero of $\opdef 1$ and $0_2$ for the trivial product of the deformation
operad which is the $0$ element in $\opdef 2(M,S)$.

Notice that neither $\op_S(M)$ nor $\op_S(M) +\op_{\mathrm{def}}(M,S)$ is a  linear operad in the sense that, although the
compositions are multilinear, the spaces for each degrees are not vector spaces but affine spaces.
On the other hand the spaces for each degrees of the deformation operad
$\op_{\mathrm{def}}(M,S)$ are vector spaces but the induced operad
compositions are not linear in general.

We may however introduce the Gerstenhaber bracket of the deformation operad
\begin{gather*}
[,]:\opdef k(M,S) \times \opdef l(M,S) \longrightarrow \opdef {k+l-1}(M,S)
\end{gather*}
defined by 
\begin{gather}\label{Gerstenhaber}
[F,G] = F \circ G - (-1)^{(k-1)(l-1)} G \circ F 
\end{gather}
where $$F \circ G = \sum_{i=1}^k (-1)^{(i-l)(l-1)} 
R(F;0_1,\dots,0_1,\underbrace{G}_{i^{\textrm{th}}},0_1,\dots,0_1).$$
This bracket is not bilinear. 
An important fact concerning this bracket is that,
\begin{gather*}
\frac{1}{2}[\gamma,\gamma]=R(\gamma;\gamma,0_1)-R(\gamma;0_1,\gamma),
\end{gather*}
which means that $\gamma$ is a product in the deformation operad iff 
\begin{eqnarray}\label{prodeq}\frac12[\gamma,\gamma] & = &  0.\end{eqnarray}

Moreover, we may define an equivalent of the Hochschild differential
\begin{gather*}
d : \opdef n(M,S) \longrightarrow \opdef {n+1}(M,S),
\end{gather*}
\begin{multline}\label{hochdiff}
dF   :=     [0_2,F]  
      =   R(0_2;F,0_1)+(-1)^{n-1}R(0_2;0_1,F)-\\
         -(-1)^{n-1}\sum_{i=1}^n(-1)^{i-1}R(F;0_1,\dots,0_1,\underbrace{0_2}_{i^{th}},0_1,\dots,0_1).
\end{multline}

It turn out that $d$ is still a coboundary operator.

\begin{Proposition}
$d$ defined by equation (\ref{hochdiff}) is a coboundary operator, i.e., $d^2= 0$.
Moreover, $\gamma\in\opdef 2(M,S)$ satisfies product equation 
$\frac12 [\gamma,\gamma] = 0,$ in $\op_{\mathrm{def}}(M,S)$ iff 
\begin{eqnarray}\label{nonlinMCeq}
d\gamma +\gamma(\gamma,S_0^1,)-\gamma(S_0^1,\gamma) & =&  0.
\end{eqnarray}
\end{Proposition}

\begin{proof}
Using equation (\ref{operadeq}) we obtain $d$ in terms of the endomorphism
compositions 
\begin{multline*}
dF  =  S_0^2(F,S_0^1)+(-1)^{n-1}S_0^2(S_0^1,F)-\\
        -(-1)^{n-1}\sum_{i=1}^n(-1)^{i-1}F(S_0^1,\dots,\underbrace{S_0^2}_{i^{th}},\dots,S_0^1).
\end{multline*}
The result follows directly from the linearity of the compositions in the endomorphism operad.
Using again equation (\ref{operadeq}) we get,
\begin{multline*}
\frac12[\gamma,\gamma] 
 =  R(\gamma;\gamma,0_1)-R(\gamma;0_1,\gamma)
 =  S_0^2(\gamma,S_0^1)+\\
 +\gamma(S_0^2,S_0^1)+\gamma(\gamma,S_0^1) -S_0^2(S_0^1,\gamma)-\gamma(S_0^1,S_0^2)-\gamma(S_0^1,\gamma),
\end{multline*}
which gives equation (\ref{nonlinMCeq}).
\end{proof}

A formal deformation $S_\epsilon$ of $S$ is a formal power series
$$S_{\epsilon} =  \epsilon S_1 + \epsilon^2 S_2 + \dots\in 
\op^n_{\textrm{form}}(M,S) 
:= \epsilon\op_{\textrm{def}}^n(M,S) \otimes k [[\epsilon]],\quad n\in\N_*,$$ 
where $\epsilon$ is a formal parameter and $\op_{\textrm{def}}(M,S)$ is
a deformation operad for $S$, 
such that $S+S_\epsilon$ is a product in $\op_S(M) + \op_{\mathrm{form}}(M,S)$.

Equivalently, one may say that $S_\epsilon$ must satisfy 
$$[S_{\epsilon},S_{\epsilon}]=0,$$
or, thanks to equation (\ref{nonlinMCeq}) that the $S_i$'s satisfy  at each order $n\in \N_*$ the following recursive
equation:
\begin{eqnarray}\label{recursive} dS_n +H_n(S_{n-1},\dots,S_1) & =  & 0,\end{eqnarray}
where $$H_n(S_{n-1},\dots,S_1) = \sum_{n= i+j} S_i(S_j,S_0^1)-S_i(S_0^1,S_i).$$

\subsubsection*{The Kontsevich product deformation}

Consider the category of real vector spaces. In this category we take 
the real vector space $M = C^{\infty}(\R^d)$ of smooth functions on $\R^d$. 
The endomorphism operad of $C^{\infty}(\R^d)$ is  
$$\op^n(M) = \Big\{\textrm{ $n$-multilinear maps from  $C^{\infty}(\R^d)^{\otimes n}$ to $C^\infty(\R^d)$}\Big\}. $$
The usual product of functions induces a product in $\op(M)$, namely
$$S_0^2(F,G)(f_1,\dots,f_k,g_1,\dots,g_l) = F(f_1,\dots,f_k)G(g_1,\dots,g_l),$$
for $F\in \op^k(M)$ and $G\in\op^l(M)$.

The induced product operad is $$\op_S^n(M) =\Big\{S_0^n\Big\},$$
where $$S_0^n(f_1,\dots,f_n) = f_1f_2\dots f_n.$$
As deformation operad, we take 
$$\opdef n(M,S) :=\Big\{\textrm{ $n$-multidifferential operators on $C^\infty(\R^d)$} \Big\}.$$
The induced coboundary operator on $\op_{\mathrm{def}}(M,S)$ is the Hochschild coboundary operator,
\begin{multline*}
dF(f_1,\dots,f_n)  =  F(f_1,\dots,f_n)f_{n+1}+(-1)^{n-1}f_1F(f_2,\dots,f_{n+1})-\\
                      -(-1)^{n-1}\sum_{i=1}^n(-1)^{(i-1)}F(f_1,\dots,f_{i-1},f_if_{i+1},f_{i+2},\dots,f_{n+1}).
\end{multline*}
and the product equation 
$$d\gamma +\gamma(\gamma,S_0^1,)-\gamma(S_0^1,\gamma)  =  0,$$
is nothing but the usual Maurer-Cartan equation. 

Kontsevich in \cite{kontsevich1997} shows that there exits a formal deformation
$$S \in \op_S^2(M)+ \epsilon \op^2_{def}(M)[[\epsilon]]$$ of $S_0^2$.
He provides the explicit formula for this deformation 
\begin{gather*}
S= S_0^2 + \sum_{n = 1}^\infty \epsilon^n \sum_{\Gamma \in G_{n,2}}
W_{\Gamma} B_{\Gamma}, 
\end{gather*}
where the $G_{n,2}$ are the Kontsevich graphs of type $(n,2)$, $W_\Gamma$ their associated
weight and $B_\Gamma$ their associated bidifferential operator (
and \cite{kontsevich1997} for more precisions).

\subsection{Monoidal structure of $\mathcal{SYM}$}

Let us recall that the extended symplectic ``category'' $\mathcal{SYM}$ is given by 
\begin{eqnarray*}
\obj & = & \Big\{\textrm{symplectic manifolds}\Big\}\\
\morph(M,N) & = & \Big\{L\subset \overline M\times N:L\textrm{ is Lagrangian}\Big\},
\end{eqnarray*}
where $\overline M$ denotes the symplectic manifold $M$ with opposite
symplectic structure $-\omega$.
The identity morphism of $\morph(M,M)$ is the diagonal $$\id_M := \Delta_M=\Big\{(m,m)\subset\overline M\times M\Big\}.$$
The composition of two morphisms $L\subset\morph(M,N)$ and $\tilde L\subset\morph(N,P)$ is
given by the composition of canonical relations,
$$\tilde L\circ L := \pi_{M\times P}\Big(
(L\times \tilde L)\cap (M\times\Delta_N\times P)
\Big)\subset \overline M\times P.$$
Everything works fine except the fact that the composition $\tilde L\circ L$ may fail to be
a Lagrangian submanifold of $\overline M\times P$. It is always the case when $L\times \tilde L$
intersects $M\times\Delta_N\times P$ cleanly (see \cite{dherin2005} for more precisions).

Let us pretend for a while that $\mathcal{SYM}$ is a true category or, better, that we have
selected special symplectic manifolds and special arrows between them such that the composition
is always well-defined.

We define the tensor product between two objects $M$ and $N$ of $\mathcal{SYM}$ as
the Cartesian product $$M\otimes N:= M\times N,$$ and the tensor product between morphisms as
\begin{eqnarray*}
L_1\otimes L_2 & := & \Big\{ (m,a,n,b):(m,n)\in L_1\\
&& \textrm{ and } (a,b)\in L_2 \Big\}\subset \morph(M\otimes A,N\otimes B),
\end{eqnarray*}
for $L_1\in \morph(M,N)$ and $L_2\in\morph(A,B)$. 

The neutral object is $\{*\}$, the one-point symplectic manifold.
The following proposition tells us that $\mathcal{SYM}$ would be a monoidal
category if it were a true category.

\begin{Proposition}
The following statements hold:
\begin{enumerate}
\item Consider $L_1\in\morph(M,A)$, $L_2\in\morph(N,B)$, $L_3\in\morph(A,X)$
and $L_4\in\morph(B,Y)$. Then
we have the following equality of sets $$(L_3\otimes L_4)\circ(L_1\otimes L_2) = (L_3\circ L_1)\otimes (L_4\circ L_2).$$
\item $\id_M\otimes \id_N = \id_{M\otimes N}$ for any object $M$ and $N$.
\item $(M\otimes A)\otimes X = M\otimes (A\otimes X)$ for any objects $M$, $A$ and $X$
\item $(L_1\otimes L_2)\otimes L_3 = L_1\otimes (L_2\otimes L_3)$ for any arrows $L_1\in\morph(M,A)$, 
$L_2\in\morph(N,B)$ and $L_3\in\morph(P,C)$.
\item $\{*\}\otimes A \simeq A\simeq A\otimes \{*\}$ for all object $A$ and 
$\id_{\{*\}}\otimes L\simeq L\simeq L\otimes \id_{\{*\}}$ for all arrows $L$, where $A\simeq B$ means that the two sets $A$ and $B$ are
in bijection.
\end{enumerate}
\end{Proposition}

\begin{proof}
(1) \begin{eqnarray*}
I & = & (L_3\otimes L_4)\circ(L_1\otimes L_2) \\
& = & \pi\Big( \big((L_1\otimes L_2)\times(L_3\otimes L_4)\big)\cap (\overline{N\times M}\times\Delta_{A\times B}\times X\times Y)\Big)\\
& = & \Big\{(m,n,\tilde x,\tilde y): \exists (a,b)\in A\times B\textrm{ s.t. }\quad  (m,n,a,b)\in L_1\otimes L_2\textrm{ and }\\
&   & (a,b,x,y)\in L_3\otimes L_4\Big\}\\ 
& = & \Big\{ (m,n,\tilde x,\tilde y):\exists a\in A,\quad (m,a)\in L_1\textrm{ and }(a,x)\in L_3\\
&   & \quad \exists b\in B,\quad (n,b)\in L_2\textrm{ and }(b,y)\in L_4 \Big\}\\
& = & (L_3\circ L_1)\otimes (L_4\circ L_2)
\end{eqnarray*}

(2) $\Delta_M\otimes\Delta_N = \Big\{(m,n,m,n):m\in M\textrm{ and } n\in N\Big\} = \Delta_{M\otimes N}.$

(3) The associativity between objects is trivial. 

(4) For morphisms, we have,

\begin{eqnarray*}
L_1\otimes L_2 & = & \Big\{(m,n,a,b):(m,a)\in L_1\textrm{ and }(n,b)\in L_2\Big\}\\
(L_1\otimes L_2)\otimes L_3 & = & \Big\{(m,n,p,a,b,c):(m,a)\in L_1,(n,b)\in L_2,\\ 
         &  & (p,c)\in L_3\Big\}
\end{eqnarray*}
and,
\begin{eqnarray*}
L_2\otimes L_3 & = & \Big\{(n,p,b,c):(n,b)\in L_2\textrm{ and }(p,c)\in L_3\Big\}\\
L_1\otimes (L_2\otimes L_3) & = & \Big\{(m,n,p,a,b,c):(m,a)\in L_1,(n,b)\in L_2, \\
  && (p,c)\in L_3\Big\}.
\end{eqnarray*}

(5) is trivial.
\end{proof}

\subsection{Lagrangian operads}

If $\mathcal{SYM}$ were a true category, we could consider the endomorphism 
operad of a 
symplectic manifold $M$. However, we may be able to restrict to a subset of Lagrangian submanifolds
$\op_{\operatorname{rest}}^n(M) \subset \op^n(M)$ for each $n\geq 0$
such that the composition 
$$L_n(L_{k_1},\dots, L_{k_n}) : = L_n\circ(L_{k_1}\otimes\dots\otimes L_{k_n}),$$
yields always a Lagrangian submanifold in 
$\op_{\operatorname{rest}}^{k_1+\dots+k_n}(M)$ for every $L_n \in 
\op_{\operatorname{rest}}^n(M)$ and
 $L_{k_i} \in \op_{\operatorname{rest}}^{k_i}(M)$, $i=1,\dots,n$. 
For instance, there is alway the trivial choice 
$$\op_{\mathrm{rest}}^1(M) = \Big\{\Delta_M\Big\},\quad \op_{\mathrm{rest}}^n(M)=\emptyset,\quad n\neq 1.$$

In this way, we may get a true operad $\op_{\operatorname{rest}}(M)$. 

The next natural question to ask is the following.

\vspace{0.5cm}
\textbf{\underline{Question:}
What is a product in a Lagrangian operad over M?}
\vspace{0.5cm}

As a first hint, take the situation where the symplectic manifold is a symplectic groupoid $G$.
In this case, we may generate an operad from the multiplication space
$G^m\in \op^2(G)$ and the base $G^{(0)}\in \op^0(G)$, the identity being the
diagonal $\Delta_G\in \op^1(G)$. Remark that $G^m$ is a product in this operad, i.e.,
that $G^m(G^m,\Delta_G) = G^m(\Delta_G, G^m)$.
Notice that the inverse of the symplectic groupoid does not play any role in this construction.

We will answer this question completely for the case were the symplectic manifold is
$T^*\R^d$ and will try to develop a deformation theory for the product in this case.

\subsubsection*{Local cotangent Lagrangian operads}

Remember that $\cotR d$ has always a structure of a symplectic groupoid over $\R^d$: the
trivial one.  The multiplication space is given in this case by 
$$\Delta_2 =\Big\{(p_1,x),(p_2,x),(p_1+p_2,x):p_1,p_2\in \R^{d*},\quad x\in \R^d\Big\}.$$
The base is $$\Delta_0 = \Big\{(0,x):x\in \R^d\Big\}.$$
If we set further 
$$\Delta_n := \Big\{(p_1,x),\dots,(p_n,x),(p_1+\dots+p_n,x):p_i\in \R^{d*},x\in \R^d\Big\},$$
it it immediate to see that the operad generated by $\Delta_0$ and $\Delta_2$ is exactly
$$\opprod n(T^*\R^d) = \Big\{\Delta_n\Big\},$$
and that $\Delta_2$ is a product in it. 

Following \cite{Cattaneo2002}, we will call this operad the \textbf{cotangent Lagrangian operad}
over $T^*\R^d$. It is the exact analog of the product operad in a monoidal category, the only
difference is that there is no true endomorphism operad to embed $\op_\Delta(\cotR d)$ into.
The idea now is to enlarge the cotangent Lagrangian operad, i.e., by considering Lagrangian
submanifolds close enough to $\Delta_n$ for each $n\in \N$ in order to have still an operad.

Notice at this point that the $\Delta_n$'s are given by generating functions. Namely, we may
identify $\overline{(T^*\R^d)^n}\times T^*\R^d$ with $T^*B_n$, where $B_n := (\R^{d*})^n\times \R^d$.
Then, 
$$
\Delta_n = \left\{
\left(\left(p_1,\frac{\partial S_0^n}{\partial p_1}(z)\right)
,\dots,
\left(p_n,\frac{\partial S_0^n}{\partial p_n}(z)\right),
\left(\frac{\partial S_0^n}{\partial x}(z),x \right)\right)
:z=(p_1,p_2,x)\in B_n
\right\}
$$
where $S_0^n$ is the function on $B_n$ defined
by 
\footnote{In the sequels, we will use the shorter notation
$(p_1+\dots+p_n)x$ instead of $\sum_{i=1}^d(p_1^i+\dots+p_n^i) x_i$.}
$$S_0^n(p_1,\dots,p_n,x) = \sum_{i=1}^d(p_1^i+\dots+p_n^i) x_i. $$
The cotangent Lagrangian operad may then be identified with
$$\opprod n =\Big\{S_0^n\Big\},\quad \opprod 0 = \Big\{0\Big\}.$$
In order to define a deformation operad for $S$, a natural idea would be to consider Lagrangian submanifolds whose generating functions are of 
the form $$F = S_0^n +\tilde F,$$ where $\tilde F\in C^\infty(B_n)$. The Lagrangian submanifold associated
to $F$ is $$L_F := \graph dF.$$
As such, the idea does not work in general.
In fact, we have to consider generating functions only defined in some neighborhood. Let us be more
precise.

We introduce the following notation, 
$$B_n^0 = \{0\}\times \R^d\subset B_n,$$ 
$V(B_n^0)$ will stand for the set of all neighborhoods of $B_n^0$ in $B_n$.

\begin{Definition}
We define $\locfuncRd n$ to be the space of germs at $B_n^0$ of smooth
functions $\tilde F$ (defined on an open neighborhood $U_{\tilde F}\subset B_n$
of $B_n^0$) which satisfy $\tilde F(0,x) = 0$ and $\nabla_p \tilde F(0,x) = 0$. 
Note that the composition will always be understood in terms of composition of germs.
\end{Definition}

\begin{Proposition}\label{prop:comp}
Let be $F\in\opprod n+\oploc n$ and $G_i\in\opprod {k_i} + \oploc {k_i}$ for $i=1,\dots,n$.
Consider the function $\phi$ defined by the formula
\begin{gather}
\Phi(p_G,x_F) = G_1\cup\dots\cup G_n(p_G,x_G) + F(p_F,x_F) - x_G p_F\label{phi-func}
\end{gather}
\begin{eqnarray*}
p_F & = & \nabla_x G_1\cup\dots\cup G_n(p_G,x_G)\label{eq1},\\
x_G & = & \nabla_p F(p_F,x_F)\label{eq2}, 
\end{eqnarray*}
where $$G_1\cup\dots\cup G_n(p_G,x_G):= G_1(p_{G_1},x_{G_1})+ \dots + G_n(p_{G_n},x_{G_n})$$ and
$p_G = (p_{G_1},\dots,p_{G_n}),$ $p_{G_i} \in (\R^{d*})^{k_i}$, $x_{G_i}\in \R^d$
and $(p_{G_i},x_{G_i})\in U_{G_i}$, for $i=1,\dots,n$.

Then, $$\phi \in \opprod{k_1+\dots+k_n}+\oploc{k_1+\dots+k_n},\quad \textrm{and}\quad L_\phi  = L_F(L_{G_1},\dots,L_{G_n}).$$
In other words, $\opprod{} +\oploc{}$ together with the product
$$\phi = F(G_1,\dots,G_n)$$
is an operad. 

Moreover, the induced operad structure on $\oploc{}$ is given by 
$$R(\tilde F;\tilde G_1,\dots,\tilde G_n) = H,$$
where $H$ is the function $H\in \oploc{k_1+\dots+k_n}$ defined by
$$H(p_G,x_F) = \tilde G(p_G,x_G) +\tilde F(p_F,x_F) -\nabla_p\tilde F(p_F,x_F)\nabla_x\tilde G(p_G,x_G),$$
\begin{eqnarray*}
p_F & = & p_F^0 + \nabla_x\tilde G(p_G,x_G),\quad p_F^0:= (p_{G_1}^\Sigma,\dots,p_{G_n}^\Sigma),\\
x_G & = & x_G^0 +\nabla_p\tilde F(p_F,x_F),\quad x_G^0 := (x_F,\dots,x_F).
\end{eqnarray*}

\end{Proposition}

\begin{Remark}[Saddle point formula]
Formula \eqref{phi-func} for $\Phi$ can be interpreted in terms of saddle point evaluation
for $\hbar\to0$ of the following integral:
\begin{multline*}
\int \ee^{\frac i\hbar
\left[F(p^1,\dots,p^k,x)+\sum_{i=1}^k \left(G_i(\pi^{i1},\dots,\pi^{il_i},y_i)
-p^i\cdot y_i\right)\right]}\;
\prod_{i=1}^k\frac{\dd^np^i\,\dd^ny_i}{(2\pi\hbar)^n} =\\
= \ee^{\frac\ii\hbar\Phi(\pi^{11},\dots,\pi^{1l_1},\pi^{21},\dots,\pi^{2l_2},\dots\dots,
\pi^{k1},\dots,\pi^{kl_k},x)}\,(C+O(\hbar)),
\end{multline*}
where $C$ is some constant.
\end{Remark}

\begin{proof}[Proof of Prop. \ref{prop:comp}]
To simplify the computations, we identify $(T^*\R^d)^n$ with $T^*(\R^{dn})$ and
$(T^*\R^d)^{k_i}$ with $T^*(\R^{dk_i})$.
With this identifications the graphs of $F$ and $G_i$, $i=1,\dots,n$  may be written as
\begin{eqnarray*}
L_F & = & \Big \{ \Big(\big(p_F, \nabla_p F(p_F,x_F)\big), \big(\nabla_xF(p_F,x_F), x_F \big) \Big):\\
    & & (p_F,x_F)\in U_F\Big\}\subset T^*(\R^{dn})\times T^*\R^d,\\
L_{G_i} & = & \Big \{ \Big(\big(p_{G_i}, \nabla_p {G_i}(p_{G_i},x_{G_i})\big), \big(\nabla_x{G_i}(p_{G_i},x_{G_i}), x_{G_i}\big) \Big):\\
        &   & (p_{G_i},x_{G_i})\in U_{G_i}\Big\}\subset T^*(\R^{dk_i})\times T^*\R^d,\\
\end{eqnarray*}
where $U_F\in V(B_n^0)$ and $U_{G_i}\in V(B_{k_i}^0)$ for $i=1,\dots,n$. 

Consider now the composition,
$$L_F(L_{G_1},\dots,L_{G_n}) = L_F \circ (L_{G_1} \otimes \dots \otimes L_{G_n}).$$
First of all, observe that,
\begin{eqnarray*}
L_G  & := &  L_{G_1} \otimes \dots \otimes L_{G_n}\\
     &  = &  \Big\{\Big(\big(p_G, \nabla_p G(p_G,x_G)\big),\big(\nabla_x G(p_G,x_G), x_G \big)\Big)
  :(p_{G_i},x_{G_i})\in U_{G_i}\Big\}\\
 L_G & \subset & T^*(\R^{d(k_1+\dots+k_n)})\times T^*(\R^{dn}).
\end{eqnarray*}
Thus,
\begin{eqnarray*}
  L_F \circ L_G & = & \pi \Big( (L_G \times L_F) \cap ( T^*\R^{d(k_1+\dots+k_n)} \times \Delta_{T^*\R^{dn}} \times T^*\R^d)\Big)\\
                & = &  \Big \{ \Big (\big( p_G, \nabla_p G(p_G,x_G)\big) ,\big( \nabla_x F(p_F,x_F),x_F\big)\Big):\\
                &   & :\quad x_G = \nabla_p F(p_F,x_F),\quad p_F = \nabla_x G(p_G,x_G), \quad(p_G,x_F)\in \tilde U \Big \}\\
L_F \circ L_G & \subset &  T^*(\R^{d(k_1+\dots+k_n)})\times T^*\R^d,  
\end{eqnarray*}

where $\tilde U$ is the subset of $(p_G,x_F)\in B_{k_1+\dots+k_n}$ such that the system,
\begin{eqnarray*}
p_F & = & \nabla_x G(p_G,x_G),\\
x_G & = & \nabla_p F(p_F,x_F),
\end{eqnarray*}
has a unique solution $(p_F,x_G)$ and such that $(p_{G_i},x_{G_i})\in U_{G_i}$, $i=1,\dots,n$, and
$(p_F,x_F)\in U_F$.  Let us check that $\tilde U$ always exists and 
is a neighborhood of $B_{k_1+\dots+k_n}^0$. To begin with, observe  that
for any $(0,x_F)\in B_n^0$ this system has the unique solution $(0,\nabla_pF(0,x_F))$. 
Set now,
\[
H(p_G,x_F,p_F,x_G) = \left(\begin{array}c 
p_F-\nabla_x G(p_G,x_G)\\
x_F-\nabla_p F(p_F,x_F)
\end{array}\right).
\]
Thanks to the fact that $G(0,x) = \sum_{i=1}^n G_i(0,x) = 0 $ we get that the Jacobi matrix
\[
D_{p_F,x_G}H\big((0,x_f,0,\nabla_pF(0,x_F)\big) = \left(\begin{array}{cc}
\id & 0\\
-\nabla_p\nabla_pF(0,x_F) & \id\end{array}\right)
\]
is invertible.

Thus, the implicit function theorem gives us the desired neighborhood $\tilde U$ of
$B_{k_1+\dots+k_n}^0$. 

Now, take $\phi$ as defined in (\ref{phi-func}). The previous considerations tell us that
$\phi$ is exactly defined on $\tilde U$.
Let us compute its graphs, 
$$L_\phi = \Big \{ \Big(\big(p_G, \nabla_p \Phi(p_G,x_F)\big), \big(\nabla_x \Phi(p_G,x_F),x_F \big)
:(p_G,x_F)\in \tilde U\Big\}.$$
We have that
\begin{multline*}
\nabla_p\phi(p_G,x_F)  =  \nabla_p G(p_G,x_G)+\nabla_xG(p_G,x_G)\frac{dx_G}{dp}+\\
	                 +\nabla_pF(p_F,x_F)\frac{dp_F}{dp}-p_F\frac{dx_G}{dp}-\frac{dp_F}{dp}x_G
		 =  \nabla_pG(p_G,x_G).
\end{multline*}
Similarly, $\nabla_x\phi(p_G,x_F) = \nabla_xF(p_F,x_F)$.
Thus, $L_\phi = L_F\circ L_G$.

At last, let us check that $\phi\in\oploc {k_1+\dots+k_n}$. First of all, remember that
\begin{eqnarray*}
F(p_F,x_F) & = & p_F^\Sigma x_F+\tilde F(p_F,x_F) \\
G(p_G,x_F) & = & \sum_{i=1}^n p_{G_i}^\Sigma x_{G_i} +\tilde G(p_G,x_G).
\end{eqnarray*}
Thus, we obtain immediately that
$$\phi(p_G,x_F) = p_G^\Sigma x_F+H(p_G,x_F),$$
where $H$ is a function only defined on $\tilde U$ by the equations,
$$H(p_G,x_F) = \tilde G(p_G,x_G) +\tilde F(p_F,x_F) -\nabla_p\tilde F(p_F,x_F)\nabla_x\tilde G(p_G,x_G),$$
\begin{eqnarray*}
p_F & = & p_F^0 + \nabla_x\tilde G(p_G,x_G),\quad p_F^0:= (p_{G_1}^\Sigma,\dots,p_{G_n}^\Sigma),\\
x_G & = & x_G^0 +\nabla_p\tilde F(p_F,x_F),\quad x_G^0 := (x_F,\dots,x_F).
\end{eqnarray*}
But now, if we set $p_G = 0$ then $p_F = 0$, $x_G = x_G^0+\nabla_p\tilde F(0,x_F)$
and $H(0,x_F) = 0$.
Similarly, one easily checks that $\nabla_p H(0,x_F) = 0 $.

\end{proof}

We will call the operad $\opprod{}+\oploc{}$
\textbf{local cotangent Lagrangian operad} over $\cotR d$ or for short
the local Lagrangian operad when no ambiguities arise. 
The induced operad $\oploc{}$ will be called the \textbf{local deformation operad of
$\opprod{}$}.

\subsubsection*{Associative products in the local deformation operad}

We say that a generating function $S\in C^\infty(B_2)$ satisfies the
\textbf{Symplectic Groupoid Associativity equation} 
if for a point $(p_1,p_2,p_3,x)\in B_3$ sufficiently close to $B_3^0$ the following implicit system for $\bar x,\bar p,\tilde x$ and $\tilde p$,
$$\bar x =\nabla_{p_1}S(\bar p ,p_3,x) ,\quad \bar p =\nabla_{x}S(p_1 ,p_2,\bar x), $$
$$\tilde x =\nabla_{p_2}S(p_1,\tilde p,x) ,\quad \tilde p =\nabla_{x}S(p_2 ,p_3,\tilde x),$$
has a unique solution and if the following additional equation holds
$$S(p_1,p_2,\bar x) + S(\bar p,p_3,x)-\bar x\bar p = S(p_2,p_3,\tilde x)+S(p_1,\tilde p,x)-\tilde x\tilde p .$$

If $S$ also satisfies the \textbf{Symplectic Groupoid Structure conditions}, i.e., if 
$$S(p,0,x)=S(0,p,x)=px \qquad\textrm{ and }\qquad S(p,-p,x) =0$$
then $S$ generates a Poisson structure  $$\alpha(x) = 2\big(\nabla_{p_k^1}\nabla_{p_l^2}S(0,0,x)\big)_{k,l=1}^d$$ on $\R^d$ 
together with a local symplectic groupoid integrating it, whose structure maps are given by
\[
\begin{array}{cccc}
\epsilon(x) & = & (0,x)  &\textbf{unit map}\\
i(p,x)      & = & (-p,x) &\textbf{inverse map}\\
s(p,x)      & = &   \nabla_{p_2}S(p,0,x) & \textbf{source map}\\
t(p,x)      & = &   \nabla_{p_1}S(0,p,x) & \textbf{target map}.
\end{array}
\]

In this case, we call $S$ a generating function of the 
Poisson structure $\alpha$ or a generating function of the local symplectic groupoid. 
See \cite{CDF2005}, \cite{dherin2004} and \cite{dherin2005} for
proofs and explanations about generating functions of Poisson structures.

The following Proposition explains what is a product in the
local cotangent Lagrangian operad.

\begin{Proposition}
 $\tilde S\in\oploc 2$ is  a product in $\oploc{}$ iff $S = S_0^2+\tilde S$ satisfies the Symplectic Groupoid
 Associativity equation.  \end{Proposition}

\begin{proof}
We know that $\tilde S$ is a product in $\oploc{}$ iff $S = S_0^2+\tilde S$ is a  product in
$\opprod{}+ \oploc{}$, i.e., iff $S(S,I) = S(I,S)$.
Let us compute.
\begin{eqnarray*}
 S(S,I)(p_1,p_2,p_3,x) & = & S\cup I(p_1,p_2,p_3,\bar x_1,\bar x_2)+S(\bar p_1,\bar p_2,x) -\bar x_1\bar p_1-\bar p_2\bar x_2\\
  & = & S(p_1,p_2,\bar x_1)+p_3\bar x_2 + S(\bar p_1,\bar p_2,x)-\bar p_1\bar x_1-\bar p_2\bar x_2,
\end{eqnarray*}
with
\begin{eqnarray*}
\bar p_1 & = & \nabla_{x_1}S\cup I(p_1,p_2,p_3,\bar x_1,\bar x_2) = \nabla_xG(\bar x)\\
\bar p_2 & = & \nabla_{x_2}S\cup I(p_1,p_2,p_3,\bar x_1,\bar x_2) = p_3\\
\bar x_1 & = & \nabla_{p_1}S(\bar p_1,\bar p_2,x)\\
\bar x_2 & = & \nabla_{p_2}S(\bar p_1,\bar p_2,x).
\end{eqnarray*}
Then we get 
$$S(S,I) = S(p_1,p_2,\bar x)+S(\bar p,p_3,x)-\bar p\bar x,$$
\begin{eqnarray*}
\bar x & = & \nabla_{p_1}S(\bar p,p_3,x)\\
\bar p & = & \nabla_xS(p_1,p_2,\bar x).
\end{eqnarray*}
Similarly, we get 
$$S(I,S) = S(p_2,p_3,\tilde x)+S(p_1,\tilde p,x)-\tilde p\tilde x,$$
\begin{eqnarray*}
\tilde x & = & \nabla_{p_2}S(p_1,\tilde p,x)\\
\tilde p & = & \nabla_xS(p_2,p_3,\tilde x).
\end{eqnarray*}
Hence, $\tilde S\in\oploc 2(T^*\R^d)$ is a product iff $S_0^2+\tilde S$ satisfies the
SGA equation. 
\end{proof}

At this point, we may still introduce the Gerstenhaber bracket as in $(\ref{Gerstenhaber})$
and the product equation in terms of the bracket would still
be $\frac12[\tilde S,\tilde S] = 0$. We may also
still write a formula for the coboundary operator. But, as this time the compositions in
$\opprod{}+\oploc{}$ are not multilinear, we cannot develop the expression $\frac12[\tilde S,\tilde S]$
in terms of the coboundary operator. Nevertheless, in Section \ref{FormalOperad}, we will develop the bracket
with help of Taylor's expansion and recover a form very close to  Equations (\ref{recursive}) in the
additive category case.

\subsubsection*{Equivalence of associative products}

To each $F\in \opprod1 +\oploc 1$,
we may associate a symplectomorphism $\psi_F$ 
which is  defined only on a neighborhood $U_F$ of $B_1^0$ in $T^*\R^d$ and which fixes $B_0^1$.
The composition of two such $\psi_G$ and $\psi_F$, which may always be defined
on a possibly smaller neighborhood $\tilde U\subset U_G$ of $B_1^0$ , is exactly 
$\psi_{F(G)}$  where $F(G)$ is the composition of $F$ by $G$ in the local Lagrangian operad.

We denote by $F^{-1} \in \opprod1+\oploc1$ the generating function of the $(\psi_F) ^{-1}$, i.e.,
the generating function such that $F(F^{-1}) = F^{-1}(F) = I$. 
Two associative products $S$ and $\tilde S$ will be called
 equivalent if $$\tilde S = F(S)(F^{-1},F^{-1})$$ for a certain $F\in\opprod1 +\oploc 1$.
It is clear that if $S\in \opprod1 + \oploc1$ is an associative product, then
$\tilde S $ also is. The following questions naturally arises.

\vspace{0.5cm}
{\bf \underline{Questions:} If $S$ generates a local symplectic groupoid, does $\tilde S$ also
generate one? Are this two local groupoids isomorphic?
}
\vspace{0.5cm}

In fact, two equivalent associative products, which are also generating functions of local
symplectic groupoids, induce isomorphic local symplectic groupoids. The isomorphism
is given explicitly by $\psi_F$. As a consequence the induced Poisson structures on the base
are the same, i.e., 
$$\alpha(x) = \nabla_{p_1}\nabla_{p_2} S(0,0,x) = \nabla_{p_1}\nabla_{p_2} \tilde S(0,0,x).$$

The following two Propositions prove these statements.

\begin{Proposition}
Let be $F\in \opprod 1+ \oploc 1$.
The following implicit equations,
\begin{eqnarray}\label{mapeq1}
x_1 & = & \nabla_p F(p_1,x_2)\\
 \label{mapeq2} p_2 & = &\nabla_x F(p_1,x_2),
\end{eqnarray}
define a symplectomorphism $\psi_F(p_1,x_1) = (p_2,x_2)$ on a neighborhood $U_F$
of $B_1^0 =\left\{(0,x):x\in\R^d\right\}$ in $T^*\R^d$ which fixes $B_1^0$ and 
which is close to the identity in the sense that $F(p,x)= px+\tilde F(p,x)$
induces the identity if $\tilde F = 0$.
Consider now $\psi_F$ and $\psi_G$ defined respectively on $U_F$ and $U_G$ for $F,G\in\opprod1 + \oploc 1$.
Then we have that $\psi_G\circ\psi_F = \psi_{F(G)}$ on  $U_{F(G)}$. 

\end{Proposition}

\begin{proof}
(1) Let us check that the system (\ref{mapeq1}) and (\ref{mapeq2}) generates a diffeomorphism around
$B_1^0$. Namely one verifies that $(\bar p_1,\bar x_1,\bar p_2,\bar x_2):= (0, \nabla_p F(0,x_2),0,x_2)$
is a solution of the system. Set now 
$$H(p_1,x_1,p_2,x_2) := \left(\begin{array}{c}x_1-\nabla_pF(p_1,x_2)\\
p_2-\nabla_xF(p_1,x_2)\end{array}\right).$$
As 
$$D_{p_1,x_1}H(\bar p_1,\bar x_2,\bar p_2,\bar x_2) =
\left(\begin{array}{cc}-\nabla_p\nabla_pF(0,\bar x_2)& \id\\
\nabla_x\nabla_pF(0,\bar x_2) & 0
\end{array}\right)$$
and 
$$D_{p_2,x_2}H(\bar p_1,\bar x_2,\bar p_2,\bar x_2) =
\left(\begin{array}{cc}0 & \nabla_x\nabla_pF(0,\bar x_2)\\
\id & 0 \end{array}\right),$$
the implicit function theorem gives us the result. Let us call $\tilde U$ the
neighborhood of $B_1^0$ where $\psi_F$ is defined.

(2) We check now that $\psi_F$ is symplectic. From equations (\ref{mapeq1}) and (\ref{mapeq2})
we get the relation $$ \frac{\partial p_l^2}{\partial p_k^1} 
= \frac{\partial x_k^1}{\partial x_l^2},$$ which directly implies that $d\psi_F \operatorname J (d\psi_F)^* = \operatorname J$
where $$\operatorname J = \left(\begin{array}{cc}0 & \id\\
-\id & 0 \end{array}\right).$$

(3) Let us see that $\psi_F(0,x) = (0,x)$. We have already noticed that 
$(0,\nabla_pF(0,x_2),0,x_2)$ is a solution of the system $(\ref{mapeq1})$ and $(\ref{mapeq2})$.
But $F(p,x) = px + \tilde F(p,x)$ with $\nabla_p\tilde F(0,p) = 0$ and then $\nabla_x\nabla_p F(0,x_2) = x_2$.

(4) Clearly $F(p,x) = px$ generates the identity.

(5) Recall that 
\begin{eqnarray*}
L_G & = & \Big \{ \Big(p_1, \nabla_p G(p_1,x_2), \nabla_xG(p_1,x_2), x_2 \Big): (p_1,x_2)\in U_G\Big\},\\
L_F & = & \Big \{ \Big(p_2, \nabla_p F(p_2,x_3), \nabla_xF(p_2,x_3), x_3 \Big):(p_2,x_3)\in U_F\Big\}.
\end{eqnarray*}
Thus, $L_G = \graph \psi_G$ and $L_F = \graph \psi_F$. The composition of these two canonical
relations yields that $L_F\circ L_G =\graph \psi_F\circ \psi_G$. 
 On the other hand, 
$L_F\circ L_G = L_{F(G)} = \graph \psi_{F(G)}$. 
Taking care on the domain of definitions, we have that 
$\psi_F\circ \psi_G = \psi_{F(G)}$ on $U_{F(G)}$.
\end{proof}

\begin{Proposition}
Let $S\in\opprod2 + \oploc2$ be a generating function of a symplectic groupoid, i.e., 
$$S(S,I) = S(I,S),\quad S(p,0,x) = S(0,p,x) = px\quad \textrm{ and }\quad S(p,-p,x) = 0.$$ Let $F \in \opprod1 + \oploc1$ such that 
$F(-p,x) = -F(p,x)$. Then,
$$\tilde S := (F(S))(F^{-1},F^{-1})$$ 
is also a generating function of a symplectic groupoid. 
The subset of odd function in $p$ forms a subgroup of $\opprod1 + \oploc 1$. 
Moreover, $\psi_F$ is a  groupoid isomorphism between the local symplectic groupoid generated
by $S$ and the one generated by $\tilde S$. As a consequence $S$ and $\tilde S$ induce the same Poisson structure
on the base.
\end{Proposition}

\begin{proof}

To simplify the notation, we set $G=F^{-1}$.
A straightforward computation gives that
\begin{multline*}
F(S)(G,G)(p_1,p_2,x) = 
S(\bar p,\tilde p,\dot x)+\\ 
+ F(\dot p,x) + G(p_1,\bar x) + G(p_2,\tilde x) -\bar p\bar x-\tilde p\tilde x-\dot x\dot p
\end{multline*}

\begin{align*}\label{myeq}
\dot x & = \nabla_p F(\dot p,x)                & \bar x & = \nabla_{p_1}S(\bar p,\tilde p,\dot x) & \tilde x  & = \nabla_{p_2}S(\bar p,\tilde x,\dot x) \\ 
\dot p & = \nabla_x S(\bar p,\tilde p, \dot x) & \bar p & =  \nabla_xG(p_1,\bar x)            &\tilde p   & = \nabla_x G(p_2,\tilde x)
\end{align*}

(1) Setting $p_1 = p$ and $p_2 = 0$, we have immediately
$$
F(S)(G,G)(p,0,x) =  G(p,\dot x) + F(\dot p,x) -\dot x\dot p
$$
with $\dot x = \nabla_{p} F(\dot p,x)$ and $\dot p = \nabla_x G(p,\dot x).$
We recognize then that
$$F(S)(G,G)(p,0,x) = F(G)(p,x) = I(p,x) = px.$$
The case $p_1 = 0$ and $p_2 = p$ is analog.

(2) One reads directly from the equation 
$$ px =  F^{-1}(p,\dot x) + F(\dot p,x) -\dot x\dot p$$
where $\dot x = \nabla_p F(\dot p,x)$ and $\dot p = \nabla_x F^{-1}(p,\dot x)$, 
that if $F$ is odd in $p$ then is also $F^{-1}$ and reciprocally.
Similarly, we check directly from the composition formula that 
$F(G)$ is odd in $p$ if $F$ and $G$ both are. Thus, the odd functions form
a subgroup of $\opprod1 + \oploc1$.

(3) Suppose now that $p_1 = p$ and $p_2 = -p$. $G$ odd in $p$ implies that 
$\bar p = -\tilde p$. As $S(p,-p,0) = 0$, we get immediately that $\tilde x =\bar x$
and $\dot p =0$ which in turns implies that $\dot x = x$.
Putting everything together, we get that $(F(S))(G,G)(p,-p,x) =0$

(4) Let us prove now that $\psi_F$ is also a groupoid isomorphism. Consider
the multiplication space of the symplectic groupoid generated by an generating function $S$, i.e, 
$$G^{(m)}( S) = \left\{
(p_1,\nabla_{p_1} S), (p_2,\nabla_{p_2} S),(\nabla_x S,x):\quad
p_1,p_2\in (\R^d)^*,\quad x \in \R^d
\right\},
$$
where the partial derivative are evaluated in $(p_1,p_2,x)$.

We have to show that $(\psi_F\times\psi_F\times\psi_F) \left(G^{(m)}(S)\right) = G^{(m)}(\tilde S)$.

A straightforward computation gives that
\begin{align*}
\nabla_{p_1}\tilde S(p_1,p_2,x) & = \nabla_{p}G(p_1,\bar x) \\
\nabla_{p_2}\tilde S(p_1,p_2,x) & = \nabla_{p}G(p_2,\tilde x)\\ 
\nabla_x \tilde S(p_1,p_2,x) & = \nabla_xF(\dot p, x).
\end{align*}

From this, we check immediately that
\begin{align*}
\psi_G\left(\left(p_1,\nabla_{p_1}\tilde S(p_1,p_2,x)\right)\right) = \left(\bar p,\nabla_{p_1}S(\bar p, \tilde p, \dot x)\right)\\
\psi_G\left(\left(p_2,\nabla_{p_2}\tilde S(p_1,p_2,x)\right)\right) = \left(\tilde p,\nabla_{p_2}S(\bar p, \tilde p, \dot x)\right)\\
\psi_F\left(\left(\nabla_x S(\bar p, \tilde p, \dot x),\dot x\right)\right) = \left(\nabla_x \tilde S(p_1,p_2,x), x\right)
\end{align*}
which ends the proof.
\end{proof}

\begin{Remark}
Suppose that $S$ is a generating function of a local symplectic groupoid. Let $F\in\opprod1 + \oploc1$
act on $S$, i.e., $\tilde S =(F(S))(F^{-1},F^{-1})$. Then, the condition $S(p,0,x) = S(0,p,x) = px$ is
preserved by any $F\in\opprod1+\oploc1$. However, the condition $S(p,-p,0)$ is only preserved by the odd
$F$s. Observe now that we have imposed the inverse map to be $i(p,x) = (-p,x)$. This implies that
$$\left(
\left(-p_2,\nabla_{p_2}S(p_2,p_1,x)\right),
\left(-p_1,\nabla_{p_1}S(p_2,p_1,x)\right),
\left(-\nabla_x S(p_2,p_1,x),x\right)
\right)\in G^{(m)}(S),$$
and thus, that $S(p_1,p_2,x) = - S(-p_1,-p_2,x)$. 
From this last equation, we get that $S$ must satisfy $S(p,-p,x) = 0$ and that
the induced local symplectic groupoid is a {\textbf symmetric} one, i.e., $t(p,x) = s(-p,x)$. 
Thus, odd transformations map symmetric groupoids to symmetric groupoids. However, they are
not the only ones.

\end{Remark}

\section{The combinatorics}\label{combinatorics}

In this Section, we present some tools which will allow us to write down at all orders the perturbative
version of the composition, Equation (\ref{phi-func}), in the local cotangent operad. All
these compositions have essentially the same form. We will first give an abstract version of the equations describing
the compositions, then we will introduce some trees which will help us to keep track of the terms involved
in the computations and, at last, we will perform the expansion in the general case.

The tools and methods presented here are essentially the same as those used in the Runge--Kutta theory 
of ODEs to determine the order conditions of a particular numeric method. We follows approximatively the notations of
\cite{GeomInt}.

\subsection{The equation}\label{equation}

Let $F:\R^{n*}\rightarrow \R^n$ and $G:\R^n\rightarrow \R^{n*}$ be two smooth functions.
Consider the point $\phi\in\R$ defined by 
\begin{eqnarray}\label{phi-eq}\phi & :=  & G(\bar x)+F(\bar p)-\bar p\bar x,\end{eqnarray}
where $\bar x$ and $\bar p$ are defined by the implicit equations,
\begin{eqnarray}\label{p-eq}
\bar p & = & \nabla_xG(\bar x)\\ \label{x-eq}
\bar x & = & \nabla_pF(\bar p).
\end{eqnarray}
Without any assumptions on $F$ and $G$, equations (\ref{p-eq}) and (\ref{x-eq}) may not 
have a solution at all or the solution may be not unique. Hence, the value $\phi$ is not
always defined.
However, if we assume that $F$ and $G$ are formal power series of the form
$$G(x)  =   p_0x+\sum_{i=1}^\infty \epsilon^i G^{(i)}(x),\quad\textrm{ and }\quad
F(p)  =   x_0p+\sum_{i=1}^\infty \epsilon^i F^{(i)}(p),$$

equations (\ref{p-eq}) and (\ref{x-eq}) become,
$$\bar p  =  p_0+\sum_{i=1}^n\epsilon^i\nabla_x G^{(i)}(\bar x),\quad\textrm{ and }\quad 
  \bar x  =  x_0+\sum_{i=1}^n\epsilon^i\nabla_p F^{(i)}(\bar p),$$
which are always recursively uniquely solvable.

Let us compute the first terms of $\bar p$, $\bar x$ and $\phi$ to get
a feeling of what is happening:
\begin{eqnarray*}
\bar p & = & p_0+\epsilon\nabla_x G^{(1)}(x_0)+\epsilon^2 \nabla_x^{(2)}G^{(1)}(x_0)\nabla_pF^{(1)}(p_0)+\cdots \\
\bar x & = & x_0+\epsilon\nabla_p F^{(1)}(x_0)+\epsilon^2 \nabla_p^{(2)}F^{(1)}(x_0)\nabla_xG^{(1)}(x_0)+\cdots \\
\phi   & = &p_0x_0+\epsilon(G^{(1)}(x_0)+F^{(1)}(p_0))+\epsilon^2 2\nabla_p F^{(1)}(p_0)\nabla_x G^{(1)}(x_0)+\cdots
\end{eqnarray*}

As we continue the expansion, the terms get more and more involved and, very soon, expressions as such become untractable.
One common strategy in physics as in numeric analysis is to introduce some graphs to keep track of the fast growing terms. Let us present
these graphs. We mainly take our inspiration from the book \cite{GeomInt}.

\subsection{The trees} \label{Cayleytree}

\begin{Definition}
-
        \begin{enumerate}
                \item{A \textbf{graph} $t$ is given by a set of vertices
		$V_t = \{1,\dots,n)$ and
                        a set of edges $E_t$ which is a set of pairs of elements of $V_t$. We denote the number
                        of vertices by $|t|$. An \textbf{isomorphism} between two graphs $t$ and $t'$ having
                        the same number of vertices is 
                        a permutation $\sigma\in S_{|t|}$ such that $\{\sigma(v),\sigma(w)\}\in E_{t'}$ if
                        $\{v,w\}\in E_{t}$. Two graphs are called \textbf{equivalent} if
                        there is an isomorphism between them. The \textbf{symmetries} of a graph are
                        the automorphisms of the graph. We denote the group of symmetries of a graph $t$ by $sym(t)$.}

                \item{A \textbf{tree} is a graph which has no cycles. Isomorphisms and symmetries are
                        defined the same way as for graphs }

                \item{A \textbf{rooted tree} is a tree  with one distinguished vertex called root. An 
                        \textbf{isomorphism} of rooted trees is an isomorphism of graphs which sends the
                        root to the root. Symmetries and equivalence are defined correspondingly.}

                \item{A \textbf{bipartite graph} is a graph $t$ together with a map $\omega:V_t\rightarrow
                        \{\circ,\bullet\}$ such that $\omega(v)\neq\omega(w)$ if $\{v,w\}\in E_t$. An isomorphism
                        of bipartite trees is an isomorphism of graphs which respects the coloring, i.e., $\omega(\sigma
                        (v))=\omega(v)$.}
		\item A \textbf{weighted graph} is a graph $t$ together with a weight map $L:V_t\rightarrow \N\backslash\{0\}$.
			An isomorphism of weighted graph is an isomorphism of graph $\sigma$ which respects the weights, i.e.,
			$\sigma(L(v)) = L(\sigma(v))$. We denote by $\|t\|$ the sum of the weights on all vertices of $t$.
        \end{enumerate}
\end{Definition}

The following table summarizes some notations we will use in the sequel.

\vspace{0.3cm}
\begin{tabular}{|c|l|}
\hline
$T$             & the set of bipartite trees\\
$RT$            & the set of rooted bipartite trees\\
$RT_\circ$      & the set of elements of $RT$ with white root\\
$RT_\bullet$    & the set of elements of $RT$ with black root\\
\hline
\end{tabular}
\vspace{0.3cm}

We will give the name {\bf Cayley trees} to trees in $T$.

We denote by $[A]$ the  set of equivalence classes of graphs in $A$ (ex: $[RT]$). They are called
\textbf{topological} ``$A$'' trees. Moreover, we denote by $A_\infty$ the weighted version of graphs in $A$. Notice
that we will use the notation $[A]_\infty$ instead of the more correct $[A_\infty]$. 
\par

The elements of $[RT]_\infty$ can be described recursively as follows:
\par
\begin{enumerate}
        \item $\circ_i,\bullet_j\in [RT]_\infty$ where $i = L(\circ_i)$ and $j= L(\bullet_j)$

        \item if $t_1,\dots,t_m\in [RT_\circ]_\infty$, then the tree $[t_1,\dots,t_m]_{\bullet_i}\in [RT]_\infty$ where
                 $[t_1,\dots,t_m]_{\bullet_i}$ is  defined by 
                connecting the roots of $t_1,\dots,t_m$ with the weighted vertex $\bullet_i$ and declaring that $\bullet_i$ is
                the new root. And the same if we interchange $\circ$ and $\bullet$. 
\end{enumerate}

Now, let us describe in terms of trees the expressions arising in the expansions of Subsection \ref{equation}.

\begin{Definition}
Given two collections of functions $F=\{F_i)_{i=1}^\infty$ and $G=\{G_j)_{j=1}^\infty$, where 
$F_i:\R^{n*}\rightarrow \R^d$ and $G_j:\R^n\rightarrow \R^{d*}$ are smooth functions, we may associate to
any rooted tree $t\in[RT]_\infty$ a vector field on $T^*\R^d$, $DC_t(F,G)\in \vect(T^*\R^d)$, called
the \textbf{elementary differential} and a function on $T^*\R^d$, $C_t(F,G)\in C^\infty (T^*\R^d)$, called
the \textbf{elementary function}.

\begin{enumerate}
\item The \textbf{elementary differential} $DC_t(F,G)$ is recursively defined as follows:
	\begin{enumerate}
	\item $DC_{\circ_i}(F,G)(p,x) = \nabla_x G^{(i)}(x)$ , $DC_{\bullet_j}(F,G)(p,x) = \nabla_p F^{(j)}(p)$
	\item $DC_{t}(F,G) = \nabla_x^{(m+1)}G^{(i)}(DC_{t_1}(F,G),\dots,DC_{t_m}(F,G))$ if $t= [t_1,\dots,t_m]_{\circ_i}$
	\item $DC_{t}(F,G) = \nabla_p^{(m+1)} F^{(j)}(DC_{t_1}(F,G),\dots,DC_{t_m}(F,G))$ if $t= [t_1,\dots,t_m]_{\bullet_j}$.
	\end{enumerate}

\item The \textbf{elementary function} $C_t(F,G)$, are recursively defined as follows: 
	\begin{enumerate}
	\item $C_{\circ_i}(F,G)(p,x) =  G^{(i)}(x)$ ,\quad $C_{\bullet_j}(F,G)(p,x) =  F^{(j)}(p)$
	\item $C_{t}(F,G) = \nabla_x^{(m)}G^{(i)}(DC_{t_1}(F,G),\dots,DC_{t_m}(F,G))$ if $t= [t_1,\dots,t_m]_{\circ_i}$.
	\item $C_{t}(F,G) = \nabla_p^{(m)} F^{(j)}(DC_{t_1}(F,G),\dots,DC_{t_m}(F,G))$ if $t= [t_1,\dots,t_m]_{\bullet_j}$.
	\end{enumerate}
\end{enumerate}
The notation $\nabla_x^{(m)}$ (resp. $\nabla_p^{(m)}$) stands for the $m^{th}$ derivative in the direction $x$ 
(resp. $p$).
\end{Definition}

Some examples are given in the following table:

\vspace{0.5cm}
\begin{tabular}{|c|c|c|}
\hline
Diagram & Elementary Differential                           &  Elementary Function \\
&&\\
\hline
\rtree{RT-1} & $\nabla_x^{(2)}G^{(i)}\nabla_p F^{(j)} $ & $\nabla_x G^{(i)}\nabla_p F^{(j)}$ \\ 
&&\\
\rtree{RT-2} &$ \nabla_p^{(3)}F^{(i)}(\nabla_x G^{(j)}, \nabla_x G^{(k)}) $& $\nabla_p^{(2)}F^{(i)}(\nabla_x G^{(j)}, \nabla_x G^{(k)})$  \\ 
&&\\
&&\\
\rtree{RT-3} &$ \nabla_x^{(3)}G^{(i)}(\nabla_p F^{(j)},\nabla_p^{(2)}F^{(k)}\nabla_x G^{(l)} ) $ & $ \nabla_x^{(2)}G^{(i)}(\nabla_p F^{(j)},\nabla_p^{(2)}F^{(k)}\nabla_x G^{(l)} ) $   \\
&&\\
\hline
\end{tabular}
\vspace{0.5cm}

Remark that for elementary functions it is not important which vertex is the root. 
This is not the case for elementary differentials. 

\begin{Definition}[\textbf{Butcher product}]
Let $u = [u_1,\dots,u_k], v = [v_1,\dots,v_l]\in [RT]$ (resp. $\in [RT]_\infty$).
We denote by
\begin{eqnarray*}
u\circ v &=& [u_1,\dots,u_k,v]\\
v\circ u &=& [v_1,\dots,v_l,u]
\end{eqnarray*}
the Butcher product. We have not written the obvious conditions
on the $u_i$ and $v_i$ so that the product remains bipartite (resp. weighted bipartite).
\end{Definition}

\begin{Definition}[\textbf{Equivalence relation on (weighted) rooted topological trees}]
Recall that an equivalence relation on a set $A$ is a special subset $R$ of $A\times A$. The equivalence
relations on $A$ are moreover ordered by inclusion. It makes then sense to consider the minimal
equivalence on $A$ containing a certain subset $U\subset A$. 

We consider here the minimal equivalence relation on $[RT]$ (resp. on $[RT]_\infty)$) such that
$u\circ v \sim  v\circ u$. 
\end{Definition}

\textbf{Properties of this relation:}
\par
It is clear that
\begin{enumerate}
\item Two topological rooted trees are equivalent if it is possible to pass from one to the
other by changing the root. More precisely:  $t,t'\in[RT]_{(\infty)}$, $t\sim t'$ iff there exists a 
representative $(E,V,r)$ of $t$ and a representative $(E',V',r')$ of $t'$ and a vertex $r''\in V$ such
that $(E,V,r'')$ and $(E',V',r')$ are isomorphic (weighted) rooted trees.
\item The quotient of $[RT]_{(\infty)}$ by this equivalence relation is exactly $[T]_{(\infty)}$.
\item It follows immediately from the definition that $C_{t}(F,G) = C_{t'}(F,G)$ if $t\sim t'$ for $i=1,2$.
\end{enumerate}

Then, it makes sense to define the elementary functions on bipartite trees.

At last, we introduce some important functions on trees: the symmetry coefficients.
\begin{Definition}
Let $t = [t_1,\dots,t_m]\in [RT]_{\infty}$. Consider the list $\tilde t_1,\dots,\tilde t_k$ of all non isomorphic
trees appearing in $t_1,\dots,t_m$. Define $\mu_i$ as the number of time the tree $\tilde t_i$ appears in
$t_1,\dots,t_m$. Then we introduce the \textbf{symmetry coefficient} $\sigma(t)$  of 
$t$ by the following recursive definition:
\par
$$\sigma(t) = \mu_1!\mu_2!\dots\sigma(\tilde t_1)\dots\sigma(\tilde t_k)$$
and initial condition $\sigma(\circ_i)=\sigma(\bullet_j)= 1$.
\end{Definition}

It is clear that $\sigma(t)$ is the number of symmetries for each representative of $t$ (i.e. $\sigma(t) = |Sym(t')|$ for
all $t'\in t$).

\subsection{The expansion} We give now a power series expansion for equation (\ref{phi-eq}).

\begin{Proposition}\label{prop:expansion}
Suppose we are given the following formal power series in $\epsilon$,
$$G(x)  =   p_0x+\sum_{i=1}^\infty \epsilon^i G^{(i)}(x),\quad \textrm{ and }\quad
F(p)  =   x_0p+\sum_{j=1}^\infty \epsilon^j F^{(j)}(p),$$

where $G^{(i)}:\R^n\rightarrow \R^{n*}$ and $F^{(j)}:\R^{n*}\rightarrow \R^n$ are smooth
functions for $i,j> 0$.

Define $\phi(p_0,x_0)\in \R[[\epsilon]]$ as 
$$\phi(p_0,x_0)  :=   G(\bar x)+F(\bar p)-\bar p\bar x,$$
where the formal power series $\bar x(\epsilon)$ and $\bar p(\epsilon)$ are uniquely determined by the implicit equations,
$$\bar p  =  p_0+\sum_{i=1}^\infty \epsilon^i\nabla_x G^{(i)}(\bar x),\quad\textrm{ and }\quad
\bar x  =  x_0+\sum_{j=1}^\infty \epsilon^j\nabla_p F^{(j)}(\bar p).$$
Then, we have that 
$$\phi(p_0,x_0) = p_0x_0 +\sum_{t\in T_\infty} \frac{\epsilon^{\|t\|}}{|t|!}C_t(F,G)(p_0,x_0).$$
\end{Proposition}

The proof of Proposition \ref{prop:expansion} is broken into several lemmas.

The method used is essentially the same as in numerical analysis when one wants
to express the Taylor series of the numerical flow of  a Runge--Kutta method.
Namely, the defining equations for $\bar p(\epsilon)$ and $\bar x(\epsilon)$ have a form very close to 
the partioned implicit Euler method(see \cite{GeomInt}).

\begin{Lemma}
 There exist unique formal power series for $\bar x(\epsilon)$ and for $\bar p(\epsilon)$  which
 satisfy equations (\ref{p-eq}) and (\ref{x-eq}). They are given by
 \begin{eqnarray}
 \label{x-sol}\bar{x}(\epsilon) &=&x_0+\sum_{t\in [RT_{\bullet}]_\infty}\frac{\epsilon^{\|t\|}}{\sigma(t)}DC_{t}(F,G),\\
 \label{p-sol}\bar{p}(\epsilon) &=&p_0+\sum_{t\in [RT_{\circ}]_\infty}\frac{\epsilon^{\|t\|}}{\sigma(t)}DC_{t}(F,G).
 \end{eqnarray}
\end{Lemma}

\begin{proof}
Uniqueness is trivial. Let us check that we have the right formal series. We only check 
equation (\ref{x-sol}). The other computation is similar.
\begin{equation*}
\begin{split}
\bar{x}(\epsilon) & =  x_0+\sum_{i\geq 1}\epsilon^i\nabla_pF^{(i)}(\bar p) \\
                  & =  x_0+\sum_{i\geq 1}\epsilon^i\sum_{m\geq 0}\frac{1}{m!}
		       \nabla_p^{(m+1)}F^{(i)}\bigg(\sum_{t\in [RT_{\circ}]_\infty}\frac{\epsilon^{\|t\|}}{\sigma(t)}DC_{t}(F,G),\dots\\
                  & \qquad \dots,\sum_{t\in [RT_{\circ}]_\infty}\frac{\epsilon^{\|t\|}}{\sigma(t)}DC_{t}(F,G)\bigg)\\
                  & =  x_0+\sum_{i\geq 1}\sum_{m\geq 0}\sum_{t_1\in [RT_\circ]_\infty}\dots\sum_{t_m\in [RT_\circ]_\infty}\frac{\epsilon^{i+\|t_1\|+\dots+\|t_m\|}}{m!\sigma(t_1) \dots\sigma(t_m)}\times\\
                  &  \qquad\times\nabla_p^{(m+1)}F^{(i)}(DC_{t_1}(F,G),\dots,DC_{t_m}(F,G))\\
                  & =  x_0 +\sum_{i\geq 1}\sum_{m\geq 0}\sum_{t_1}\dots\sum_{t_m}\frac{\epsilon^{\|t\|}}%
		  {m!\sigma(t)}(\mu_1!\mu_2!\dots) DC_{t}(F,G),\\
		  &    \quad \textrm{with }t=[t_1,\dots,t_m]_{\bullet_i}\\
                  & =  x_0 + \sum_{t\in[RT_\bullet]_\infty}\frac{\epsilon^{\|t\|}}{\sigma(t)}  DC_{t}(F,G).
\end{split}
\end{equation*}
\end{proof}

\begin{Lemma} We have the following expansion for $\phi(p_0,x_0)$:
\begin{multline*}
\phi(p_0,x_0) = p_0x_0+ \sum_{t\in [RT]_\infty} \frac{\epsilon^{\|t\|}}{\sigma(t)} C_{t}(F,G)-\\
-\Big(\sum_{t\in [RT_\circ]_\infty} 
\frac{\epsilon^{\|t\|}}{\sigma(t)} DC_{t}(F,G) \Big) \Big(\sum_{t\in [RT_\bullet]_\infty}
\frac{\epsilon^{\|t\|}}{\sigma(t)} DC_{t}(F,G) \Big).
\end{multline*}
\end{Lemma}

\begin{proof}
We compute the different terms arising in $G(\bar x)+F(\bar p)-\bar p\bar x$ in terms of
trees.
\begin{equation*}
\begin{split}
 G(\bar x)  & =  p_0\bar x + \sum_{i\geq 1}\epsilon^i\sum_{m\geq 0}\frac{1}{m!}\nabla_x^{(m)}G^{(i)}%
            \bigg(\sum_{t\in [RT_{\bullet}]_\infty}\frac{\epsilon^{\|t\|}}{\sigma(t)}DC_{t}(F,G) ,\dots\\
            & \qquad   ,\dots,\sum_{t\in [RT_{\bullet}]_\infty}\frac{\epsilon^{\|t\|}}{\sigma(t)}DC_{t}(F,G)\bigg)\\
            & =  p_0\bar x + \sum_{i\geq 1}\sum_{m\geq 0}\sum_{t_1\in [RT_\bullet]_\infty}\dots%
	    \sum_{t_m\in [RT_\bullet]_\infty}\frac{\epsilon^{\|t\|}}{m!\sigma(t)}(\mu_1!\mu_2!\dots)\times\\
            &    \quad  \times\nabla_x^{(m)}G^{(i)}(DC_{t_1}(F,G),\dots,DC_{t_m}(F,G)),\\
	    &    \quad\qquad \textrm{with }t=[t_1,\dots,t_m]_{\bullet_i}\\ 
            & =  p_0\bar x + \sum_{t\in [RT_\circ]_\infty}\frac{\epsilon^{\|t\|}}{\sigma(t)} C_{t}(F,G)
\end{split}
\end{equation*}

By the same sort of computations we obtain,
$$F(\bar p)  =  x_0\bar p+\sum_{t\in [RT_\bullet]_\infty} \frac{\epsilon^{\|t\|}}{\sigma(t)} C_{t}(F,G).$$
Finally, we get the desired result as,
\begin{multline*}
p_0\bar x+x_0\bar p-\bar p\bar x = p_0x_0 -\\
-\Big(\sum_{t\in [RT_\circ]_\infty} 
\frac{\epsilon^{\|t\|}}{\sigma(t)} DC_{t}(F,G) \Big) \Big(\sum_{t\in [RT_\bullet]_\infty}
\frac{\epsilon^{\|t\|}}{\sigma(t)} DC_{t}(F,G) \Big).
\end{multline*}
\end{proof}

Thus, $\phi(p_0,x_0)$ is expressed as sums over topological weighted rooted bipartite trees. We would like now to regroup
the terms of the formula in the previous Lemma. To do so, we express all terms in terms of topological trees (no longer rooted).

\begin{Lemma}
Let $u\in[RT_\circ]_\infty$ and $v\in [RT_\bullet]_\infty$. Then, 
$$DC_{u}(F,G)DC_{v}(F,G) = C_{u\circ v}(F,G) = C_{v\circ u}(F,G).$$
\end{Lemma}

\begin{proof}
Suppose $u=[u_1,\dots,u_m]_{\circ_i}$, $v=[v_1,\dots,v_l]_{\bullet_j}$, then we get
\begin{eqnarray*}
A & = & DC_{u}(F,G)DC_{v}(F,G)\\
 & = & \nabla_x^{(m+1)}G^{(i)}(DC_{u_1}(F,G),\dots,DC_{u_m}(F,G)).DC_{v}(F,G)\\
               & = & \nabla_x^{(m+1)}G^{(i)}(DC_{u_1}(F,G),\dots,DC_{u_m}(F,G),DC_{v}(F,G))\\
               & = & C_{u\circ v}(F,G).\\
\end{eqnarray*}

\end{proof}

\begin{Lemma}\label{Lemma:sym}
Let $t=(V_t,E_t)\in T_\infty$. For all $v\in V_t$ let $t_v$ be the bipartite rooted tree $(V_t,E_t,v)\in RT_\infty$.
For $v\in V_t$ and $e= \{u,v\}\in E_t$ we have

\begin{eqnarray*}
 \frac{|sym(t)|}{|sym(t_v)|}  & = & |\{v'\in V_t /t_{v'}\textrm{is isomorphic to } t_v\}|\\
 \frac{|sym(t)|}{|sym(t_u)||sym(t_v)|}  & = & |\{e'\in E_t /t_{u'}\sqcup t_{v'}\textrm{is isomorphic to } t_u\sqcup t_v\}|
\end{eqnarray*}
\end{Lemma}

\begin{proof}
Consider the induced action of the symmetry group of the tree on the set of
vertices. Notice that two vertices $v$ and $w$ are in the same orbit iff
$t_v$ is isomorphic to $t_w$. Then the number of vertices of $t$ which lead to rooted
tree isomorphic to $t_v$ is exactly the cardinality of the orbit of $v$, which is exactly
$|sym(t)|$ divided by the cardinality of the isotropy subgroup which fixes $v$. But the latter
is $|sym(t_v)|$ by definition. We then get the first statement.
\par
For the second statement we have to consider the induced action on the edges and apply the
same type of argument. 
\end{proof}

\begin{Lemma}
We get
$$\phi(p_0,x_0) = p_0x_0+\sum_{t\in T_\infty}\frac{\epsilon^{\|t\|}}{|t|!}C_t(F,G).$$
\end{Lemma}

\begin{proof}
Let us perform the last computation.
\begin{equation*}
\begin{split}
\phi(p_0,x_0) & = p_0x_0+ \sum_{t\in [RT]_\infty} \frac{\epsilon^{\|t\|}}{\sigma(t)}C_{t}(F,G)-\\
             &\qquad -\sum_{u\in[RT_\circ]_\infty}\sum_{v\in[RT_\bullet]_\infty} \frac{\epsilon^{\|u\|+\|v\|}}{\sigma(u)\sigma(v)}DC_{u}(F,G)DC_{v}(F,G)\\
        & = p_0x_0+ \sum_{\bar{t}\in [T]_\infty} \epsilon^{|\bar{t}|}C_{\bar{t}}(F,G)
                \Big\{ \sum_{t\in\bar{t}} \frac{1}{|sym(t)|}-\\
		 &\qquad-\sum_{\substack{u\in [RT_\bullet]_\infty,v\in [RT_\circ]_\infty\\ u\circ v \in\bar{t}}}%
                   \frac{1}{|sym(u)||sym(v)|}   \Big\}\\
        & = p_0x_0+ \sum_{t\in T_\infty}\frac{\epsilon^{\|t\|}}{|t|!}C_{t}(F,G)\Big\{%
        \sum_{v\in V_t} \frac{|sym(t)|}{|sym(t_v)|}\frac{1}{k(t,v)}-\\
         &\qquad -\sum_{e = \{u,v)\in E_t} \frac{|sym(t)|}{|sym(t_u)||sym(t_v)|}\frac{1}{l(t,e)} \Big\}
\end{split}
\end{equation*}
where $k(t,v) = |\{v'\in V_t /t_{v'}\textrm{is isomorphic to } t_v\}|$ 
and $l(t,e) = |\{e'\in E_t /t_{u'}\sqcup t_{v'}\textrm{is isomorphic to } t_u\sqcup t_v)|$.
Using Lemma \ref{Lemma:sym} and the fact that for a tree the difference between the number of vertices
and the number of edges is equal to 1 we get the desired result.
\end{proof}

Using now the fact that $S$ is a formal power series we immediately get Proposition \ref{prop:expansion}.

\section{Deformation of a non-linear structure} \label{FormalOperad}

\subsection{The formal cotangent Lagrangian operad}

The formal cotangent Lagrangian operad on $T^*\R^d$ is the perturbative/formal version of the local cotangent operad on $T^*\R^d$. 
Recall that in the latter the product for $F \in \opprod n+\oploc n$ 
and $G_i \in \opprod n+ \oploc {k_i}$, $i= 1,\dots n$ was expressed as in Proposition \ref{prop:comp}:
$$F(G_1,\dots,G_n)(p_G,x_F)  =  G_1 \cup \dots \cup G_n (p_G,x_G) + F(p_F,x_F) - p_F \cdot x_G, $$
\begin{eqnarray*}
p_F & = & \nabla_x G_1 \cup \dots \cup G_n (p_G,x_G) , \\
x_F & = & \nabla_p F(p_F,x_F). 
\end{eqnarray*}
If we consider $p_G$ and $x_F$ as parameters in the previous equations, we have then that,
$$G(p_G, \cdot)  :  \R^{nd} \longrightarrow \R,\quad\textrm{ and } F( \cdot, x_F) :  (\R^{nd})^*  \longrightarrow \R, $$
Suppose now that the $F$ and $G_i$, $i=1,\dots,n$, are formal series of the form
\begin{eqnarray*}
F(p_F,x_F) & = & p_F^{\Sigma} \cdot x_F + \sum_{i=1}^{\infty} \epsilon^i F^{(i)}(p_F,x_F)\\
G_l(p_{G_l},x_{G_l}) & = & p_{G_l}^{\Sigma} \cdot x_{G_l} + \sum_{i=1}^{\infty}\epsilon^i G_l^{(i)}(p_{G_l},x_{G_l})
\end{eqnarray*}
where  $$p^{\Sigma} := \sum_{i=1}^n p_i\quad\textrm{ for }\quad p=(p_1,\dots,p_n)\in(\R^{dn})^*.$$ 
We may rewrite $F$ and $G$ as,
\begin{gather*}
F(p_F,x_F) = x_0^Fp_F + \sum_{i=1}^{\infty} \epsilon^i F^{(i)}(p_F,x_F)
\end{gather*}
\begin{gather*}
G(p_G,x_G) = p_0^G x_G + \sum_{i=1}^{\infty} \epsilon^i G^{(i)}(p_G,x_G)
\end{gather*}
where $x_0^F = (x_F,\dots,x_F) \in \R^{dn}$ and 
$p_0^G = (p_{G_1}^\Sigma,\dots,p_{G_n}^\Sigma) \in (\R^{nd})^*$ for  $x_G \in \R^{dn}$ and $p_F \in (\R^{dn})^*.$

Applying now Proposition \ref{prop:expansion}, we obtain for the compositions the following expansion:
\begin{multline}\label{treecomp}
F(G_1,\dots,G_n)(p_G,x_G) = p_G^{\Sigma} \cdot x_F  +\\  +\sum_{t \in T_{\infty}} \frac{\epsilon^{||t||}}{|t|!} C_t \Big(
 F(\cdot,x_F),G_1 \cup \dots \cup G_n (p_G,\cdot)\Big) (p_0^G,x_0^F).
\end{multline}

This motivates to define the formal deformation space of the cotangent Lagrangian operad $\opprod{}(\cotR d)$ as
 $$ \opform n(\cotR d,\Delta) :=   \Big\{  \sum_{i=1}^{\infty}\epsilon^i F^{(i)} : F^{(i)} \in P_i^n(T^*\R^d) \Big\},$$
where $P_i^n(T^*\R^d)$ stands for the vector space of functions $F:B_n \longrightarrow \R$ such that
\begin{enumerate}
\item $F(p,x)$ is a polynomial in the variables $p=(p_1,\dots,p_n)$,
\item $F(\mu p,x) = \mu^{i+1}F(p,x).$
\end{enumerate}

One may think of $\oploc{} +\opform{}$  as the Taylor series of functions in $\opprod{}+\oploc{}$
The compositions are given by formula (\ref{treecomp}), which also tells us that $\opprod{}+\oploc{}$
is an operad. The unit is $$I(p,x)= px,\quad I \in \opprod{} + \opform 1.$$

The induced operad structure on $\opform{}$ is then given by,
\begin{gather*}
  I \in \opform 1,\quad  I(p,x)=0,\\
\opform n = \Big\{ \sum_{i=1}^{\infty} \epsilon^i F^{(i)} : F^{(i)} \in P_i^n(T^*\R^d)\Big \}\\
 F(G_1,\dots,G_n)(p_G,x_F) = \sum_{t \in T_{\infty}}\frac{\epsilon^{||t||}}{|t|!} C_t \big(
 F, G_1 \cup \dots \cup G_n\big).
\end{gather*}
This operad will be called \textbf{the formal deformation operad} of the cotangent
Lagrangian operad $\opprod{}$.

\subsection{Product in the formal deformation operad}

Exactly as for the local deformation operad, $S_\epsilon$ is a product in $\opform{}$ iff
$S_0^2+S_\epsilon$  satisfies formally the SGA equation. Moreover, if $S_0^2+S_\epsilon$
satisfies the SGS conditions, then $S_0^2+S_\epsilon$ is the generating function of
a formal symplectic groupoid over $\R^d$.

Again, the zero of $\opform2$ is a product in $\opform{}$. We will stick to the
conventions introduced for $\oploc{}$. Namely, $0_1$ will stand for the zero of $\opform1$,
which is also the identity of the operad and
$0_2$ will stand for the zero of $\opform2$, which is the trivial product of the operad.

Thanks to the composition formula (\ref{treecomp}), we are now able to 
rewrite the product equation in $\opform{}$ as a cohomological equation, exactly
as the deformation equation of a product in an additive category. Note that the Taylor expansion
plays the same role as the linear expansion played in the additive case.

Let us define the Gerstenhaber bracket in $\opform{}$ as follows:
$$[F,G] = F\circ G-(-1)^{(k-1)(l-1)}G\circ F,$$
where $$F\circ G = \sum_{i=1}^k (-1)^{(i-l)(l-1)}F\{0_1,\dots,0_1,\underbrace{G}_{i^{th}},0_1\,\dots,0_1),$$
for $F\in \opform k$ and $G\in \opform l$.

We are now able to define a true coboundary operator.

\begin{Proposition}
Consider  $d:\opform n\rightarrow \opform{n+1}$ 
$$dF := [0_2,F].$$ Then, $d$ may be written as 
\begin{multline*}\label{Sdiff}
dF(p_1,\dots,p_{n+1})  =  F(p_1,\dots,p_n,x)+\\
 +\sum_{j=1}^n(-1)^{n+j-1}F(p_1,\dots, p_j+p_{j+1},\dots, p_n,x) + (-1)^{n-1}F(p_2,\dots,p_{n+1},x)
\end{multline*}
Moreover,  $d$ is linear and $d^2 = 0$
\end{Proposition}

\begin{proof}
For more clarity, let us break our convention and write $\tilde I$ instead of $0_1$ and $\tilde S$ instead of $0_2$.
We have that  $[\tilde S,F] = \tilde S\circ F-(-1)^{n-1}F\circ\tilde S$. As $\tilde S = 0$, only the 
trees $\circ_i$ and $\bullet_j$ will contribute to the product. Then we have,
\begin{eqnarray*}
I_1 & = & \tilde S\circ F(p_1,\dots,p_{n+1},x)  \\
&=& \sum_{i\geq 1}\epsilon^i \Bigg(C_{\bullet_i}\bigg(\tilde S(\cdot,x),F\cup \tilde I(p,\cdot)\bigg)\Big((\sum_1^n p_l,p_{n+1}),(x,x)\Big)+  \\
&&  +(-1)^{n-1}C_{\bullet_i}\bigg(\tilde S(\cdot,x),\tilde I\cup F(p,\cdot)\bigg)\Big((p_1,\sum_2^{n+1} p_l),(x,x)\Big)\Bigg)\\
&=&  \sum_{i\geq 1}\epsilon^i\Big(F^{(i)}(p_1,\dots,p_n,x)+ (-1)^{n-1}F^{(i)}(p_2,\dots,p_{n+1},x)\Big)
\end{eqnarray*}
and 
\begin{eqnarray*}
I_2 & = & F\circ\tilde S(p_1,\dots,p_{n+1})  \\
& = & \sum_{j=1}^n(-1)^{j-1}\sum_{i\geq 1}\epsilon^i %
C_{\circ_i}\bigg(F(\cdot,x),( \tilde I\cup\dots \tilde I \cup \underbrace{\tilde S}_{j^{th}}\cup\tilde I\dots\\
& & \dots \cup \tilde I)(p,\cdot)%
   \bigg)\Big((p_1,\dots,p_j+p_{j+1},\dots,p_{n+1}),(x,\dots,x)\Big)\\
& = & \sum_{j=1}^n (-1)^{j-1}\sum_{i\geq 1} \epsilon^i F^{(i)}(p_1,\dots,p_j+p_{j+1},\dots,p_{n+1},x),
\end{eqnarray*}
which gives the desired formula.
The check that $d^2 = 0$ is straightforward.
\end{proof}

We have then a complex $$\Big(C^\bullet= \oplus_{n\geq 0}\opform n, d\Big).$$ This complex is exactly
the Hochschild complex of (formal) multi-differential operators lifted
on the level of symbols ( see for instance \cite{CahenGutt1980}). This remark gives us the cohomology of the complex,
$$\operatorname H^n(C^\bullet) \simeq \mathcal \epsilon \mathcal V^n(\R^d)[[\epsilon]],$$
where $\mathcal V^n(\R^d)$ is the space of $n$-multi-vector fields on $\R^d$.

We come now to the question of finding a product $S_\epsilon$ in the formal deformation operad
of $\opprod{}$ This is exactly the same problem  as deforming the trivial generating function
$S_0^2$ in $\opprod{}+\opform{}$. We are thus looking for an element $ S_\epsilon\in\opform 2$ of the form
$$ S_\epsilon = \epsilon  S_1+\epsilon^2 S_2+\dots$$
such that \begin{eqnarray}\label{MaurerCartan} [ S_\epsilon, S_\epsilon] = 0. \end{eqnarray}
 Equation (\ref{MaurerCartan}) becomes, on the level of trees,
\begin{eqnarray}\label{operadSGA}\sum_{t\in T_\infty}\frac{\epsilon^{\|t\|}}{|t|!}
\Big( C_t(S_\epsilon,S_\epsilon\cup I)-C_t(S_\epsilon,I\cup S_\epsilon) \Big)
&  =  & 0.\end{eqnarray}

One sees immediately  that this equation is equivalent to the following infinite set of recursive equations,
\begin{eqnarray*}\label{cohomeq} d S_n + H_n(S_{n-1},\dots,S_1) & = & 0,\end{eqnarray*}
where 
$$H_n(S_{n-1},\dots,S_1) = \sum_{\substack{t\in T_\infty^{k,n}\\2\leq |k|\leq n}}\frac1{|t|!}
\Big( C_t(S_\epsilon,S_\epsilon\cup I)-C_t(S_\epsilon,I\cup S_\epsilon) \Big),
$$
where $T_\infty^{k,n}$ is the subset of trees in $T_\infty^{k,n}$ with $k$ vertices and
such that $\|t\| = n$. These recursive equations are the 
 exact analog of  Equations (\ref{recursive}).

\subsection{Formal symplectic groupoid generating function}
\label{proof}

We restate now the main Theorem  of \cite{CDF2005}, Theorem 1, in terms of the new structures
defined in this article.

\begin{Theorem}
For each Poisson structure $\alpha$ on $\mathbb R^d$, we have that 
$$S_\epsilon(\alpha)  = \sum_{n=1}^\infty \frac{\epsilon^n}{n!}\sum_{\Gamma\in T_{n,2}}W_\Gamma \hat B_\Gamma(\alpha)$$
is a product in the formal deformation operad $\opform{}(\cotR d,\Delta)$ of the cotangent Lagrangian operad 
$\opprod{}(\cotR d)$.  Moreover, $S_\epsilon(\alpha)$ is the unique natural product in  $\opform{}(\cotR d,\Delta)$ whose first order is $\epsilon \alpha$.
\end{Theorem}

In the above Theorem, the $T_{n,2}$ stand for the set of Kontsevich trees of type $(n,2)$, $W_\Gamma$ is the
Kontsevich weight of $\Gamma$ and $\hat B_\Gamma$ is the symbol of the bidifferential operator $B_\Gamma$ associated to $\Gamma$.
We refer the reader to \cite{CDF2005} for exact definitions of Kontsevich trees, weights, operators and naturallity.

We called $S_\epsilon(\alpha)$ the (formal) {\bf  symplectic groupoid generating function} because, as shown in \cite{CDF2005},
it generates a ``geometric object'', a (formal) symplectic groupoid over $\R^d$ associated to the Poisson structure $\alpha$
whose structure maps are explicitly given by
\[
\begin{array}{cccc}
\epsilon_\epsilon(x) & = & (0,x)  &\textbf{unit map}\\
i_\epsilon(p,x)      & = & (-p,x) &\textbf{inverse map}\\
s_\epsilon(p,x)      & = & x + \nabla_{p_2}S_\epsilon(\alpha)(p,0,x) & \textbf{source map}\\
t_\epsilon(p,x)      & = & x + \nabla_{p_1}S_\epsilon(\alpha)(0,p,x) & \textbf{target map}.
\end{array}
\]

This exhibits a strong relationship between star products and symplectic groupoids already
foreseen by Costes, Dazord, Weinstein, Karasev, Maslov and Zakrzewski in respectively \cite{CDW1987}, \cite{karasev1989} and \cite{zakrzewski1990}
. Recently and from a completely different point of view, Karabegov in \cite{karabegov2004bis} went still a step further 
by showing how to associate a kind of ``formal symplectic groupoid'' to any star product.

In \cite{dherin2004} and \cite{dherin2005}, we prove that the product $S_\epsilon(\alpha)$ has a non-zero
convergence radius provided that the Poisson structure $\alpha$ is analytic. In this case, the generated formal symplectic
groupoid is the local one. We also compared compared this local symplectic groupoid with the one constructed by Karasev 
and Maslov in \cite{karasev1989} and we proved that this two local symplectic groupoids are not only isomorphic as they should
but exactly identical.

\end{document}